\documentclass[11pt]{article}

\pdfoutput=1
\usepackage[pdftex]{color}
\usepackage{amssymb}
\usepackage{amsthm}
\usepackage{amsmath}
\usepackage{latexsym}
\usepackage{amscd}
\usepackage{graphicx}
\usepackage[pdftex, colorlinks=true, citecolor=green]{hyperref}
\usepackage{lscape}
\usepackage{multirow}

\newcommand{\ds}{\displaystyle}

\newcommand{\ben}{\begin{equation}}     
\newcommand{\eeqn}{\end{equation}}
\newcommand{\bey}{\begin{eqnarray}}
\newcommand{\eey}{\end{eqnarray}}

\newtheorem{thm}{Theorem}[section]
\newtheorem{prop}[thm]{Proposition}
\newtheorem{lemma}[thm]{Lemma}

\newtheorem{defn}[thm]{Definition}
\newtheorem{conj}[thm]{Conjecture}

\begin{document}

\begin{flushleft}
{\Large
\textbf{Template iterations of quadratic maps and hybrid Mandelbrot sets}
}
\\
\vspace{4mm}
 Anca R\v{a}dulescu$^{*,}\footnote{Assistant Professor, Department of Mathematics, State University of New York at New Paltz; New York, USA; Phone: (845) 257-3532; Email: radulesa@newpaltz.edu}$, Kelsey Butera$^1$, Brandee Williams$^1$
\\
\indent $^1$ Department of Mathematics, SUNY New Paltz, NY 12561
\\
\end{flushleft}

\vspace{3mm}

\begin{abstract}
\noindent As a particular problem within the field of non-autonomous discrete systems, we consider iterations of two quadratic maps $f_{c_0}=z^2+c_0$ and $f_{c_1}=z^2+c_1$, according to a prescribed binary sequence, which we call a \emph{template}. We study the asymptotic behavior of the critical orbits, and define the Mandelbrot set in this case as the locus for which these orbits are bounded. However, unlike in the case of single maps, this concept can be understood in several ways. For a fixed template, one may consider this locus as a subset of the parameter space in $(c_0,c_1) \in \mathbb{C}^2$; for fixed quadratic parameters, one may consider the set of templates which produce a bounded critical orbit. In this paper, we consider both situations, as well as \emph{hybrid} combinations of them,  we study basic topological properties of these sets and interpret them in light of potential applications.
\end{abstract}

\noindent {\bf Summary: Using non-autonomous iterations of discrete functions, we build a mathematical framework that can be used to study the effect of errors in copying mechanisms (such as DNA replication). In our theoretical setup -- in which one of the functions is the correct one, and the other one is the erroneous perturbation -- we consider problems that a sustainable replication system may have to solve when facing the potential for errors. We find that it is possible to tell which specific errors are more likely to affect the system's dynamics, in absence of prior knowledge of their timing. Moreover, within an optimal locus for the correct function, even a small number of errors can affect the sustainability of the system. Mathematically, our work complements broader existing results in non-autonomous dynamics with more specific detail for the case of two arbitrary iterated functions, which is a valuable context for applications.}

\section{Introduction}

Discrete dynamics of single iterated maps on the complex plane has been a rich field of studies over many decades, in particular for polynomial maps. In this context, the Julia set is defined as the boundary between initial conditions which remain asymptotically bounded and those which escape to infinity under iterations of the map. The topological and fractal properties of Julia sets have been well studied for polynomials, with major results relating the geometry of the Julia set with properties of the critical orbits~\cite{branner1992iteration,qiu2009proof,devaney2006criterion}. One of the most studied families is that of quadratic maps  in the family $f_c = z^2+c$, with $c \in \mathbb{C}$, with a history that goes back more than a century, to the work of Fatou and Julia~\cite{fatou1920equations,julia1918memoire}. For iterations of single quadratic maps, it is known that the Julia set is either connected, if the orbit of the critical point 0 is bounded, or totally disconnected, if the orbit of the critical point 0 is unbounded. The postcritically bounded parameter locus is therefore the same as the Julia set connectedness locus in the parameter complex plane, and is know as the Mandelbrot set~\cite{branner1989mandelbrot}, the topology and properties of which have been amply studied since the 1960s.

While single map iterations have been often used to represent natural phenomena, it is unlikely that natural systems evolve according to the same identical dynamics along time. A more realistic mathematical framework to model the variability and errors that appear in replication systems is that of time-dependent (random) iterations, in which the iterated map may change between steps (evolve in time). 

A broad field of active research in probability theory centers around dynamical systems produced by random iterations. For example, Kiefer et al. have been investigating theoretical ergodic properties of systems generated by iterations of functions chosen at random from a certain family, according to some probability distribution~\cite{kifer2012ergodic}. Diaconis and Freedman have developed methods for studying the steady state distribution of a Markov chain, and have given useful convergence criteria, in particular for iterates of random Lipschitz functions~\cite{diaconis1999iterated}. Other studies relate the idea of random function iterations to applications in economics (via dynamical systems subject to random shocks)~\cite{bhattacharya2007random}, or to periodically forced, monotone difference equations motivated by applications from population dynamics~\cite{cushing2002periodically,haskell2005stochastic}. In the particular case of two maps $f_0$ and $f_1$, one can consider a probability measure on the sequence space $\Sigma = \{0, 1\}^N$ which, in the simplest case, is generated as a product measure of Bernoulli probabilities $p$ and $1-p$ on the two symbols in the set $\{ 0,1 \}$. Along these lines, Bhattacharya and Rao have been studying invariant measures of Markov processes obtained by iterations of maps chosen at random from a set of two quadratic maps~\cite{bhattacharya1993random}.

A particular direction of research has been to study time-independent iterations of complex maps. Rather than focusing specifically on the structure of the probability space, the aim has been to describe general properties of the asymptotic dynamics for arbitrary sequences of maps chosen from specific families (e.g., hyperbolic polynomials of a certain degree). The focus on this aspect of non-autonomous systems was introduced in the work of Fornaess and Sibony ~\cite{fornaess1991random}, and continues in more recent work by Comerford, Stankewitz and Sumi~\cite{comerford2013preservation,comerford2006hyperbolic,sumi2010random}. The simple particular case of an iteration which alternates two distinct quadratic complex maps (thus equivalent to iterating a complex  quartic polynomial) was studied by Danca et al.~\cite{danca2009alternated}. It was shown that in this case the Julia sets can be disconnected without being totally disconnected, and that they exhibit a graphical alternation of patterns found independently within the Julia set of each of the two iterated maps.

In this paper, we continue our previous work in the spirit of this latter line of research~\cite{radulescu2015symbolic}, while incorporating a simple probability space structure on the sequence space. More specifically, we study the dynamics generated by two different complex quadratic functions, $f_{c_0}$ and $f_{c_1}$, applied according to a general binary symbolic sequence ${\bf s}$ (\emph{template}), in which the ``zero'' positions correspond to iterating the function $f_{c_0}$ and the ``one'' positions correspond to iterating the function $f_{c_1}$. We view template iterations  as a framework for replication or learning algorithms that occur in nature, with patterns that evolve in time, and which may involve occasional, random or periodic ``errors.'' While our results are generally expected to agree with existing general results in non-autonomous iterations of complex polynomials, our interest resides more specifically in understanding properties that are specific to a system based on a quadratic map pair. We study the dependence of the dynamic behavior on the two primary features of this system's hardwiring: (1) the complex parameter pair $(c_0,c_1)$ that fixes the iterated maps and (2) the structure of the template, i.e. the particular succession (timing) of the 0s and 1s that govern the iteration order.

In our previous work, we studied the parameter locus for which the orbit of the initial critical point remains bounded. More precisely, for a fixed $(c_0,c_1,{\bf s}) \in \mathbb{C}^2 \times \{ 0,1 \}^\infty$, we studied the iterated orbit:
$$\xi_0=0 \to f_{c_{s_1}}(0) \to (f_{c_{s_2}} \circ f_{c_{s_1}})(0) \to \cdots \to (f_{c_{s_n}} \circ \cdots \circ f_{c_{s_1}})(0) \to \cdots$$

\noindent When $c_0=c_1$, the orbit does not depend on ${\bf s}$ and corresponds to the traditional orbit of $\xi_0$ under the quadratic polynomial $f_{c_0} = f_{c_1}$. In this slice, the Mandelbrot set is defined as the locus of all $c_0$ for which the orbit of the origin is bounded, or equivalently for which the Julia set of $f_{c_0}$ is connected. When extending the concept to template iterations based on the parameter space $\mathbb{C}^2 \times \{ 0,1 \}^\infty$, one natural question to ask is, for example, whether the two traditional definitions are still equivalent. In this scenario, in which an ``erroneous'' function $f_{c_1}$ interferes into the iteration of a correct function $f_{c_0}$, one may also ask how far can the erroneous map $c_1$ stray from the given correct map $c_0$, so that we can be sure (with some probability) that the critical orbit $o_{\bf s}(0)$ stays bounded when a random template ${\bf s}$ prescribes where the errors will hit in the iteration process. 

To better understand such questions, our direction is to consider different slices of the full parameter space, as well as hybrid combinations; for these simpler sets, we discuss topological properties of potential importance (such as connectivity), and we explore connections with other possible definitions. The paper is organized as follows:

After discussing some basic definitions and properties in Section~\ref{intro}, we review two concepts already introduced in previous work~\cite{radulescu2015symbolic}: the \emph{fixed map} and the \emph{fixed template} Mandelbrot sets, as parameter slices in $\{0,1\}^\infty$ and in $\mathbb{C}^2$, respectively. We had previously documented some of the properties of these sets for both periodic and random templates (e.g., there are detectable differences in the Hausdorff dimension along the boundary of  Mandelbrot slices for fixed templates, between the case of periodic and that of non-periodic templates). In Section~\ref{structure}, we concentrate on fixing the complex parameter pair $(c_0,c_1)$, and further understanding the topology of the \emph{fixed-map} Mandelbrot set, defined as a Mandelbrot slice in the space $\{ 0, 1 \}^\infty$. 

In Section~\ref{hybrid_M}, we define \emph{hybrid} Mandelbrot sets to illustrate, in the form of a surface over the complex plane, the likelihood of a random template to produce a bounded critical orbit $o_{\bf s}(0)$ under iterations of a fixed map $f_{c_0}$ and a variable map $f_{c_1}$. We find that the level sets of hybrids have a less complex geometric structure than the traditional Mandelbrot set, possibly relating to a known phenomenon of cooperation between generating maps towards ``smoothing out the chaos''~\cite{sumi2010random}.

A different possibility, raised previously by Comerford~\cite{comerford2013preservation}, is to consider the Mandelbrot set to represent the parameter locus for which the iteration is post-critically bounded. This definition is more consistent with accounting for multiple critical orbits when iterating higher order polynomials. In our case, this definition means that none of the critical points $\xi_0 = 0$ introduced at any step $m \geq 1$ of the iteration can escape to infinity. In other words, all orbits:
$$0 \to f_{c_{s_m}}(0) \to \cdots \to (f_{c_{s_n}} \circ \cdots \circ f_{c_{s_m}})(0) \to \cdots$$
are bounded. In Section~\ref{multicrit} of our paper, we consider this stricter, ``multi-critical'' definition of the Mandelbrot set. We define the corresponding parameter slices and investigate the same questions as in the uni-critical definition, and compare the results. In particular, we discuss, for both definitions, the relationship of the Mandelbrot set with the connectedness locus of the Julia set for template iterations.


\section{The quadratic family and template iterations}
\label{intro}

\noindent We will be working within the complex quadratic family 
$$\{ f_c \colon \mathbb{C} \rightarrow \mathbb{C}, f_c(z) =z^2 + c,  \text{ with } c \in \mathbb{C} \}.$$ 

\noindent For fixed parameters $c_0, c_1 \in \mathbb{C}$, and a fixed binary sequence ${\bf s} = (s_n)_{n \geq 1} \in \{0,1\}^\infty$, one can define the template orbit for any $\xi_0 \in \mathbb{C}$ as the sequence $o_{\bf s}(\xi_0) = (\xi_n)_{n \geq 1}$ constructed recursively, for every $n \geq 0$, as
$$\xi_{n+1} = f_{c_{s_{n+1}}} (\xi_n).$$

\vspace{5mm}
\noindent The parameter space of our template iterated system is $\mathbb{C}^2 \times \{ 0,1 \}^\infty$, with the quadratic parameter pair $(c_0,c_1) \in \mathbb{C}^2$, and the template ${\bf s} \in \{ 0,1 \}^\infty$. 

\begin{defn}
We define the Mandelbrot set for template iterations as
$${\cal M} = \{ (c_0,c_1,{\bf s}) \in \mathbb{C}^2 \times \{ 0,1 \}^\infty , \text{ where } o_{\bf s}(0) \text{ is bounded under iterations of $f_{c_0}$ and $f_{c_1}$} \}.$$
\end{defn}

As a step towards understanding the Mandelbrot set as a whole (in the full parameter space $\mathbb{C}^2 \times \{ 0,1 \}^\infty$), one may first consider studying ``Mandelbrot'' slices, which are easier to describe and visualize. In our previous work~\cite{radulescu2015symbolic}, we defined two main types of slices:

\begin{defn}
Fix ${\bf s} \in \{ 0,1 \}^\infty$ a symbolic sequence. The corresponding \textbf{fixed-template Mandelbrot slice} is defined as
$${\cal M}_{\bf s} = \{ (c_0,c_1) \in \mathbb{C}^2, \text{ such that } (c_0,c_1,{\bf s}) \in {\cal M}  \}.$$
\end{defn}

\begin{defn}
Fix $(c_0,c_1) \in \mathbb{C}^2$. The corresponding \textbf{fixed-map Mandelbrot slice} is defined as
$${\cal M}_{c_0,c_1} = \{ {\bf s} \in \{ 0,1 \}^\infty, \text{ such that } (c_0,c_1,{\bf s}) \in {\cal M} \}.$$
\end{defn}


\noindent The space $\{ 0,1 \}^\infty$ is an ultrametric space with the metric induced by the distance
$$d({\bf s},{\bf t}) = \sum_{k=1}^\infty \frac{|s_k-t_k|}{2^k}$$

\noindent for any two templates ${\bf s}$ and ${\bf t}$. To better represent subsets of templates with certain properties, one can define an order on $\{ 0,1 \}^\infty$ by identifying each template with the binary expansion of a real number in the unit interval, using the distance function to the null space:
$$\psi \colon \{ 0,1 \}^\infty \to [0,1]$$
$$\psi({\bf s}) = \sum_{n=1}^\infty s_n 2^{-n} \text{ for all } {\bf s}=(s_n)_{n \geq 1}$$

The map $\psi$ is surjective, but not injective, since any number in $[0,1]$ has at least one binary expansion, but some numbers do have multiple expansions (e.g., $0.1\overline{0} = 0.0\overline{1}$).

Using this map, we can represent the set ${\cal M}_{c_0,c_1}$ as its image $\psi({\cal M}_{c_0,c_1}) \subset [0,1]$, which we will further illustrate and discuss in Section~\ref{structure}, where we also show that $\psi({\cal M}_{c_0,c_1})$ are measurable with the Lebesgue measure ${\cal L}$ on the unit interval. In fact, this corresponds to ${\cal M}_{c_0,c_1}$ being measurable in $\{0,1\}^\infty$, with ${\cal \mu}({\cal M}_{c_0,c_1}) = {\cal L}(\psi({\cal M}_{c_0,c_1}))$, where ${\cal \mu}$ represents the product measure induced on $\{0,1 \}^\infty$ by the equally weighted Bernoulli probability $p(0)=p(1)=1/2$ on each template trial. For the sake of concreteness, we choose to illustrate in the remainder of the section the images of the actual sets ${\cal M}_{c_0,c_1}$, as subsets of the unit interval with the Lebesgue measure.

On the $\mathbb{C}^2$ slice of the parameter space for template iterations, one can consider the function
$$b \colon \mathbb{C}^2 \to [0,1], \: \text{ given by } b(c_0,c_1) = {\cal L} \left( \psi({\cal M}_{c_0,c_1}) \right).$$

\noindent For a fixed $c_0 \in \mathbb{C}$, we define
$$b_{c_0} \colon \mathbb{C} \to [0,1], \: \text{ given by } b_{c_0}(c_1) = b(c_0,c_1).$$

\noindent One can then view the ``graph'' of each $b_{c_0}$ as a mixed, \emph{hybrid} Mandelbrot set, potentially related to phenomena of averaging out chaos, as described by Sumi in the context of non-autonomous iterations on polynomial semigroups~\cite{sumi2010random}. 

\begin{defn}
Fix $c_0 \in \mathbb{C}$. The corresponding \textbf{hybrid Mandelbrot set} is defined as:
$${\cal M}_{c_0} = \{ (c_1,b_{c_0}(c_1)) \in \mathbb{C} \times [0,1], \text{ for all } c_1 \in \mathbb{C} \}.$$
\end{defn}

\noindent While in previous work we focused on fixed-template Mandelbrot sets, in the current paper we will investigate fixed-map and hybrid sets. The main goal of this paper is to visualize the structure and conjecture properties of these sets, in particular understand their dependence on the complex parameters and on the template structure, whichever is appropriate in each case.


\subsection{Escape radius for template iterations}

A simple, yet major result in the case of iterations of single quadratic complex maps is the existence of an escape radius~\cite{branner1989mandelbrot}. More precisely:

\begin{thm}
For $\lvert c \rvert<2$, $R_e=2$ is an escape radius for $f_c(z) = z^2 +c$. In particular, for the critical orbit $(z_n)_{n \geq 1}$ (with $z_0=0$), if $\lvert z_N \rvert > 2$ for some positive integer $N$, then $\lvert z_n \rvert > 2$ for all $n \geq N$, and $\lvert z_n \rvert \to \infty$ as $n \to \infty$.
\label{escape_singlemap}
\end{thm}

\noindent More general results have been formulated by Comerford et al. for the case of non-autonomous iterations~\cite{comerford2013preservation}, as follows: Suppose we take $d \geq 2$, $K \geq 1$, $M \geq 0$ and let $(P_m)_{m=1}^\infty$ be a bounded sequence of polynomials, that is: (1) the degree $d_m$ of each $P_m$ satisfies $2 \leq d_m \leq d$; (2) the leading coefficients $a_m$ satisfy $1/K \leq a_m \leq K$; (3) all coefficients of all $P_m$ are within the disc of radius $M$. Then there exists an escape radius $R >0$ for the sequence $P_m$, that is: for all $m \geq 0$ and all $\lvert z \rvert > R$, we have $\lvert P_n \circ  \hdots \circ P_{m+2} \circ P_{m+1}(z) \rvert \to \infty$ as $n \to \infty$. Comerford et al.~\cite{comerford2013preservation} showed that one can find an escape radius $R$ which depends only on the bounds $d$, $K$ and $M$, and works for every sequence which satisfies these bounds.

The existence and the value of the escape radius will be used in this study to determine numerically which orbits are unbounded,  allowing elimination from the Mandelbrot set of the parameters for which the critical point iterates outside the escape radius. While Comerford's result permits custom computation of a tighter $R_e$ for a general sequence of bounded polynomials, when iterating two quadratic maps according to a template, it it possible to use a more convenient escape rate, reminiscent of single map iterations. Below, we show that the single map escape radius $R_e=2$ still acts as an escape radius for template iterations, under some relatively broad assumptions for the maps. We have the following:

\begin{thm} For arbitrary $(c_0,c_1,{\bf s}) \in \mathbb{C}^2 \times \{ 0,1 \}^\infty$, consider a template orbit $(\xi_n)_{n \geq 1}$. Suppose $\lvert \xi_N \rvert > R_e$, with $R_e=\max\{ 2, \lvert c_0 \rvert, \lvert c_1 \rvert \}$ for some $N$. Then $\lvert \xi_n \rvert > R_e$ for all $n \geq N$, and $\lvert \xi_n \rvert \to \infty$ as $n \to \infty$. In other words, $R_e = \max\{ 2, \lvert c_0 \rvert, \lvert c_1 \rvert \}$ is an escape radius for the template iteration of the maps $c_0$ and $c_1$.
\label{escape_lemma}
\end{thm}

\noindent {\bf Proof.} Notice that 

$$\frac{\lvert \xi_{n+1} \rvert}{\lvert \xi_n \rvert} = \frac{\lvert f_{c_{s_n}}(\xi_n) \rvert}{\lvert \xi_n \rvert} = \frac{\lvert \xi_n^2 + c_{s_n} \rvert}{\lvert \xi_n \rvert} \geq \frac{\lvert \xi_n \rvert^2 - \lvert c_{s_n} \rvert}{\lvert \xi_n \rvert} = \lvert \xi_n \rvert - \frac{\lvert c_{s_n} \rvert}{\lvert \xi_n \rvert}.$$

\noindent We will prove by induction that $\lvert \xi_n \rvert >R_e$, for $n \geq N$. We know that this is true for $n=N$, since $\lvert \xi_N \rvert >R_e$ by hypothesis. Suppose that $\lvert \xi_n \rvert >R_e$, for some $n \geq N$. We will show that $\lvert \xi_{n+1} \rvert >R_e$. Indeed, since $\lvert \xi_n \rvert > \lvert c_{s_n} \rvert$ by the induction hypothesis, we have that

$$\frac{\lvert \xi_{n+1} \rvert}{\lvert \xi_n \rvert} \geq \lvert \xi_n \rvert - \frac{\lvert c_{s_n} \rvert}{\lvert \xi_n \rvert} > \lvert \xi_n \rvert -1.$$ 

\noindent Since $\lvert \xi_n \rvert >2$, we further have:
$$\frac{\lvert \xi_{n+1} \rvert}{\lvert \xi_n \rvert} > \lvert \xi_n \rvert - 1 > 1.$$

\noindent Hence $\lvert \xi_{n+1} \rvert > \lvert \xi_n \rvert >R_e$, which concludes the induction step. In conclusion, $\lvert \xi_n \rvert > R_e$ for all $n \geq N$. 

It also follows from this that $\lvert \xi_{n+1} \rvert > \lvert \xi_n \rvert$ for all $n \geq N$, implying that the orbit increases in radius towards its upper bound. If we assume that this upper bound is a finite number $\ds l = \lim_{n \to \infty} \lvert \xi_n \rvert$, we get, taking limit on both sides of $\ds \frac{\lvert \xi_{n+1} \rvert}{\lvert \xi_n \rvert} \geq \lvert \xi_n \rvert - 1$, that the limit $l \leq 2$, which contradicts the fact that $l$ is the upper bound of a sequence of values larger than 2. It follows that the limit $\ds \lim_{n \to \infty} \lvert \xi_n \rvert = \infty$. 

\hfill $\Box$ \\

\noindent {\bf Remark.} This implies directly that, when $(c_0,c_1) \in \mathbb{C}^2$ with $\lvert c_0 \rvert, \lvert c_1 \rvert \leq 2$, $R_e=2$ is an escape radius for the template iteration corresponding to the pair of maps $(c_0,c_1)$.

\begin{prop} For fixed $c_0=0$ and any arbitrary $c_1 \in \mathbb{C}$, the template orbit $(\xi_n)_{n \geq 1}$ of $\xi_0=0$ under $(0,c_1)$ escapes once it falls outside the radius $R_e=2$. That is, if $\lvert \xi_N \rvert > 2$ for some positive integer $N$, then $\lvert \xi_n \rvert > 2$ for all $n \geq N$, and $\lvert \xi_n \rvert \to \infty$ as $n >N$.
\end{prop}

\noindent {\bf Proof.} Suppose first that $\lvert c_1 \rvert \leq 2$. From $\lvert \xi_N \rvert > 2$, it follows automatically in this case that $\lvert \xi_N \rvert > \max\{ 2, \lvert c_0 \rvert, \lvert c_1 \rvert \}$. From Theorem~\ref{escape_lemma}, it further follows that $\lvert \xi_n \rvert > 2$ for all $n \geq N$, and $\lvert \xi_n \rvert \to \infty$ as $n >N$.

Now we consider the case $\lvert c_1 \rvert > 2$. We aim to track down when the critical orbit $o_{\bf s}(0)$ leaves the disk of radius $2$. The orbit will remain zero as long as the map $f_{c_0}$ is iterated, that is throughout the first succession of zeros in the template. (If $f_{c_1}$ is never iterated, the critical orbit is trivial.) Otherwise, the first nonzero orbit entry will be $\xi_M = f_{c_1}(0) = c_1$, with $\lvert \xi_M \rvert >2$. This will be squared repeatedly over the next succession of zero template entries. If all the remaining template entries are zero, the orbit will be, for all $n \geq 0$, $\xi_{M+n} = c_1^{2^n}$, with $\lvert \xi_{M+n} \rvert > 2$. If there is at least one additional nonzero entry in the template, at position $M+K$, then $\xi_{M+K} = c_1^{2^K}$, and $\xi_{M+K+1} = c_1^{2^{K+1}}+c_1$. Notice that $\lvert \xi_{M+K+1} \rvert = \lvert c_1^{2^{K+1}}+c_1 \rvert \geq \lvert c_1 \rvert ( \lvert c_1 \rvert^{2^{K+1}-1} -1) > \lvert c_1 \rvert > 2$. Hence $\xi_{M+K+1}$ satisfies the conditions of Theorem~\ref{escape_lemma} and the orbit escapes. In either case, once the orbit escapes the disk of radius $R_e=2$, it continues to increase to $\infty$.

\hfill $\Box$

\noindent {\bf Remark.} Clearly, $R_e=2$ is not an escape radius for any arbitrary template iteration. Here is a simple construction for which $R_e=2$ fails to act as an escape radius. Fix a value of $c_0$ outside of the traditional Mandelbrot set. We choose the value of $c_1$ and the template ${\bf s}$ as follows. According to Theorem~\ref{escape_singlemap}, the orbit of zero leaves the disk of radius $R_e=2$ after a certain number $N$ of iterations of $f_{c_0}$, and never returns: $0 \to z_1 \to z_2 \to \hdots \to z_N \to$, with $\lvert z_n \rvert >2$, for all $n \geq N$. Consider now a point $z'$ in the filled Julia set of $f_{c_0}$ (hence $|z'|<2$), and choose $c_1 = z'-z_N^2$ (notice that $\lvert c_1 \rvert > |z_N|^2-|z'|>2$). Choose the template ${\bf s}$ to be such that $s_j=0$, for all $j \neq N+1$, and $s_{N+1}=1$. Then the orbit $o_{\bf s}(0)$ leaves the escape disk in $N$ iterates (since $\lvert z_N \rvert >2$), but $z_{N+1} = z_{N}^2 + c_1 = z'$, and the orbit remains henceforth bounded (all future iterations are under $f_{c_0}$). Notice that the choice that generates this counterexample is not robust, since the filled Julia set is totally disconnected for a $c_0$ outside of the Mandelbrot set; hence the choice of $z'$, and subsequently of $c_1$, is made from a set of area zero.


\section{Structure of fixed-map Mandelbrot sets}
\label{structure}

Henceforth, we will focus on parameters $c_0$ and $c_1$ in the open complex disc of radius two $\mathbb{D}(2)$. In this section, we will study the structure of fixed-map Mandelbrot sets, and how this changes when the pair $(c_0,c_1) \in \mathbb{D}(2)^2$ is varied. Since later in the paper we perform a computer-assisted numerical analysis, we will use truncated templates in order to illustrate our sets. We therefore find it useful to define finitely iterated versions of our sets (that is, Mandelbrot sets for binary templates of finite length $N$), and explore the limit behavior as $N \to \infty$.

\begin{defn}
For any integer $N \geq 1$, we call the $N$-root of the template ${\bf s} = (s_n)_{n\geq 1}$ the finite binary sequence $\langle {\bf s} \rangle_{_N} = s_1 \hdots s_N$. We will use the notation $\{ 0,1 \}^N$ for the set of all $N$-roots. 
\end{defn}

\noindent We will say that two templates ${\bf s}$ and ${\bf s}' \in \{ 0,1 \}^\infty$ have a common $N$-root if they agree up to the $N$th position: $s_j = s_j'$, for $1 \leq j \leq N$. $N$ need not be maximal in order to define a common $N$-root; if two templates have a common $N$-root, they will have a common $k$-root, for all $1 \leq k \leq N$. We will use the notation $o_{\langle {\bf s} \rangle_{_N}}(\xi_0)$ for the finite orbit of a point $\xi_0$ under the $N$-root $\langle {\bf s} \rangle_{_N}$, and $o_{\langle {\bf s} \rangle_{_N}}^j(\xi_0)$ to designate the $j$-th iterate of $\xi_0$ under the $N$-root $\langle {\bf s} \rangle_{_N}$.

\begin{defn}
For a fixed pair $(c_0,c_1) \in \mathbb{D}(2)^2$, we define the \textbf{$N$-rooted fixed-map Mandelbrot set} as:
$${\cal M}^N_{c_0,c_1} = \{ {\bf s} \in \{ 0,1 \}^\infty, \text{ with } \lvert o_{\langle {\bf s} \rangle_{_N}}^N(\xi_0) \rvert \leq 2 \}$$
\end{defn}

\vspace{2mm}
\noindent With this definition, together with the definition of the fixed-map Mandelbrot set in Section~\ref{intro}, we have the following lemma:

\begin{lemma}
For $\lvert c_0 \rvert, \lvert c_1 \rvert <2$ and for any arbitrary integer $N \geq 1$ we have that ${\cal M}_{c_0,c_1} \subseteq {\cal M}_{c_0,c_1}^{N+1} \subseteq  {\cal M}_{c_0,c_1}^N$. Moreover, $\ds \bigcap_{N=1}^\infty {\cal M}_{c_0,c_1}^N = {\cal M}_{c_0,c_1}$.
\end{lemma}

\noindent {\bf Proof.} Suppose $\langle {\bf s} \rangle_{_N}$ is an $N$-root such that $\lvert o_{\langle {\bf s}\rangle_{_N}}^N(0) \rvert > 2$ (that is, the orbit of $\xi_0=0$ escapes by the time the iteration reaches the end of the $N$-root). Then none of the templates ${\bf s} \in \{ 0,1 \}^\infty$ with the common $N$-root $\langle s \rangle_{_N}$ are in ${\cal M}_{c_0,c_1}$. Hence ${\cal M}_{c_0,c_1} \subseteq {\cal M}_{c_0,c_1}^{N+1}$. Moreover, if $\lvert o_{\langle {\bf s} \rangle_{_{N+1}}}(0) \rvert \leq 2$ for some template ${\bf s}$ then, from the escape radius condition, it follows that $\lvert o_{\langle {\bf s} \rangle_{_N}}(0) \rvert \leq 2$ as well, which implies that ${\cal M}_{c_0,c_1}^N$ are nested around ${\cal M}_{c_0,c_1}$. Subsequently, the intersection $\ds \bigcap_{N=1}^\infty {\cal M}_{c_0,c_1}^N \supseteq {\cal M}_{c_0,c_1}$.

On the other hand, notice that if  $\ds {\bf s} \in \bigcap_{N=1}^\infty {\cal M}_{c_0,c_1}^N$, so that $\lvert o_{\langle {\bf s} \rangle_{_N}}(0) \rvert \leq 2$ for all $N \geq 1$, then ${\bf s} \in {\cal M}_{c_0,c_1}$. Hence $\ds \bigcap_{N=1}^\infty {\cal M}_{c_0,c_1}^N \subseteq {\cal M}_{c_0,c_1}$. In conclusion
$$\ds \bigcap_{N=1}^\infty {\cal M}_{c_0,c_1}^N = {\cal M}_{c_0,c_1}.$$
\hfill $\Box$  \\

\begin{figure}[h!]
\begin{center}
\includegraphics[width=0.6\textwidth]{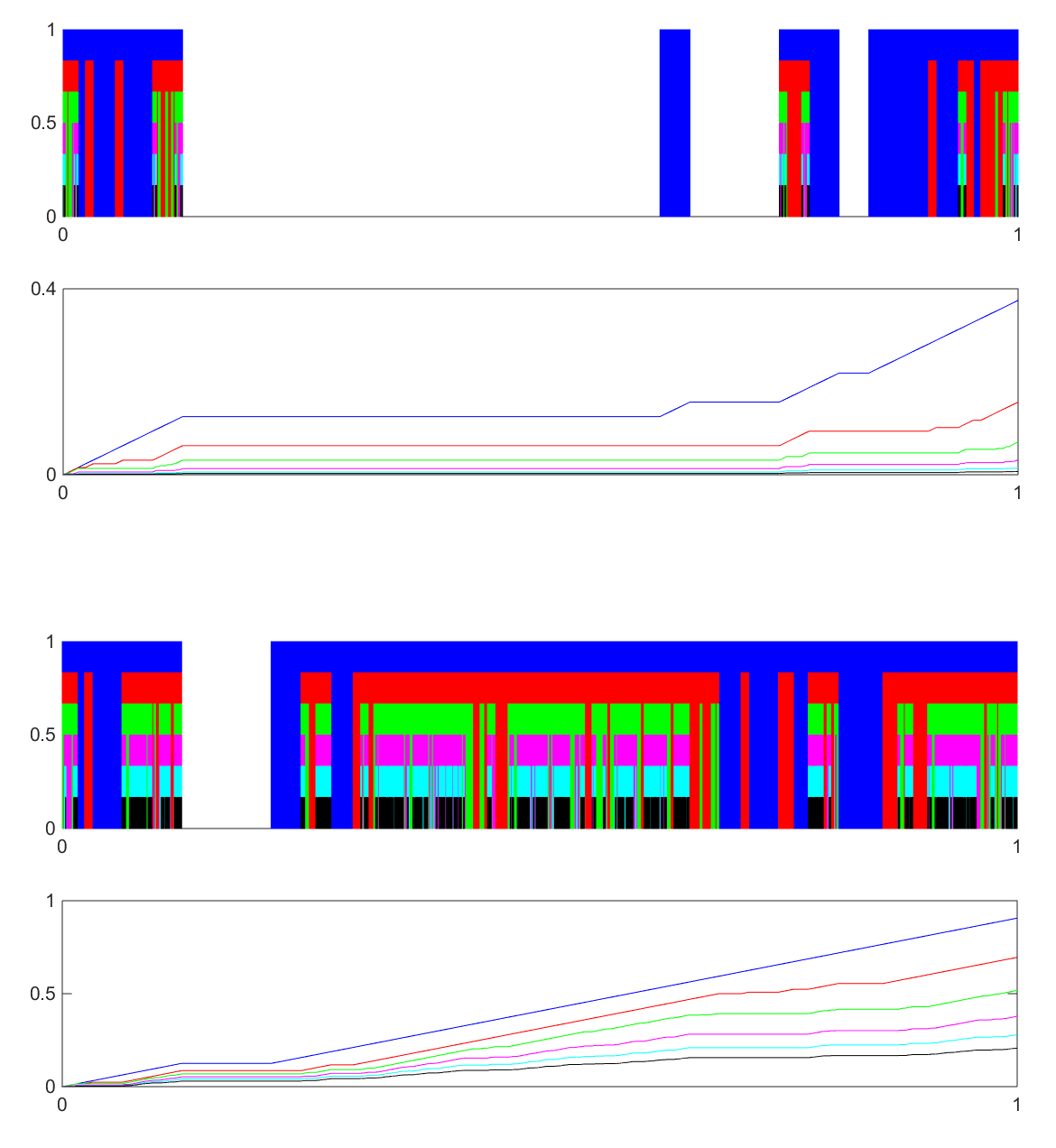}
\end{center}
\caption{\emph{\small {\bf $N$-rooted fixed-map Mandelbrot sets} for two different map pairs $(c_0,c_1)$: $c_0=-0.117-0.76i$ and $c_1=-0.62-0.62i$ ({\bf top}); $c_0=-0.75$ and $c_1=-0.117-0.856i$ ({\bf bottom}). The top subplot of each panel shows sets $\psi({\cal M}_{c_0,c_1}^N)$ for different values of $N$, in different colors: $N=5$ (blue), $N=7$ (red), $N=9$ (green), $N=11$ (pink), $N=13$ (cyan), $N=15$ (black). To make it easier to simultaneously visualize these sets within the same copy of the unit interval, the sets were represented as blocks with hight decreasing from 1 to 1/6 for increasing $N$, layering out their structure for comparison (they are nested sets). The bottom subplot represents the accumulation map for each of the sets on top, in corresponding colors.}}
\label{combs}
\end{figure}

\noindent Notice first that the intersection ${\cal M}_{c_0,c_1}$ is nonempty if either $c_0$ or $c_1$ are in the traditional Mandelbrot set, (${\cal M}_{c_0,c_1}$ contains the all zero template if $c_0$ is in ${\cal M}$, and it contains the all one template if is $c_1$ is in ${\cal M}$). ${\cal M}_{c_0,c_1}$ is empty if both $\lvert c_0 \rvert >2$ and $\lvert c_1 \rvert >2$. The ``size'' of ${\cal M}_{c_0,c_1}$ for intermediate values of $c_0$ and $c_1$ will be further discussed in Section~\ref{hybrid_M}.

Since $\psi({\cal M}_{c_0,c_1}^N)$ are Lebesgue measurable, it follows that $\psi({\cal M}_{c_0,c_1})$ is also measurable (as claimed in Section~\ref{intro}). Figure~\ref{combs} illustrates on the unit interval the nesting of $\psi({\cal M}_{c_0,c_1}^N)$ to $\psi({\cal M}_{c_0,c_1})$, for two distinct pairs $(c_0,c_1)$ of parameters in $\mathbb{C}^2$. 
For a fixed $N$, there are $2^N$ possible $N$-roots $\langle {\bf s}^0 \rangle_{_N} = 0 \hdots 00$, $\langle {\bf s}^1 \rangle_{_N} = 0 \hdots 01$, up to $\langle {\bf s}^{2^N} \rangle_{_N} = 1 \hdots 11$, whose images under $\psi$, together with the endpoint $1$, form a partition $(a^N_j)_{1 \leq j \leq 2^N+1}$ of $[0,1]$ (the dyadic fractions of level $N$): $a_j^N = \psi(\langle {\bf s}^j \rangle_{_N})$, for $j = 0,2^N$, and $a_{2^N+1}^N = 1$. It is easy to see that $\psi({\cal M}_{c_0,c_1}^N)$ is a union of intervals $[a_j^N,a_{j+1}^N]$. Similarly, $\psi({\cal M}_{c_0,c_1}^{N+1})$ is a union of intervals $[a_j^{N+1},a_{j+1}^{N+1}]$ formed by points in the partition generated by all $(N+1)$-roots. Note that $(a_j^{N+1})$ is a finer partition than $(a_j^N)$. Suppose now that an interval $[a_j^{N+1},a_{j+1}^{N+1}] \subset \psi({\cal M}_{c_0,c_1}^{N+1})$. This means that the last iterate of zero under the $(N+1)$-root $\langle {\bf s} \rangle_{_{N+1}} = \psi^{-1}(a_j^{N+1})$ is inside the escape disc. Consequently, all the previous iterates under $\langle {\bf s} \rangle_{_{N+1}} = s_1s_2 \hdots s_{_{N+1}}$ have to also be inside the escape disc. If $a_i^N = \psi(s_1,...,s_N)$, it follows that $[a_j^{N+1},a_{j+1}^{N+1}] \subset [a_i^N,a_{i+1}^N]$. Hence every interval in $\psi({\cal M}_{c_0,c_1}^{N+1})$ is part of an interval in $\psi({\cal M}_{c_0,c_1}^N)$. Figure~\ref{combs} represents the sets $\psi({\cal M}_{c_0,c_1}^N)$ in different colors for different values of $N$. In order to make the nested property more visible, we plotted them as 2-dimensional ``combs'' $\psi({\cal M}_{c_0,c_1}^N) \times [0,h]$, where the hight $h$ decreases with the value of $N$ (so that one can observe the layering of nested combs of different colors as $N$ increases). 

We define the \emph{accumulation map} for each ${\cal M}_{c_0,c_1}^N$ as follows. Consider $\{ 0,1 \}^N = \{ \langle {\bf s}^1 \rangle_{_N},\hdots \langle {\bf s}^{2^N} \rangle_{_N} \}$ the set of all $N$-roots, and their corresponding partition $(a_j^N)_{1 \leq j \leq 2^N+1}$, as defined above. Then the accumulation map is given by  \\

$\hspace{4cm} \ds \phi_{c_0,c_1}^N \colon [0,1] \to [0,1],$\\

$\hspace{4cm} \ds \phi_{c_0,c_1}^N(t)= \left \{ \begin{array}{l}
0, \text{ at } t=0\\\\
\phi(a_j^N) \text{ on } [a_j^N,a_{j+1}^N], \text{ if } \langle {\bf s}^j \rangle \notin {\cal M}_{c_0,c_1}^N \\\\
\ds \phi(a_j^N) + t-a_j^N \text{ on } [a_j^N,a_{j+1}^N], \text{ if } \langle {\bf s}^j \rangle_{_N} \in {\cal M}_{c_0,c_1}^N.
\end{array} \right. $

\vspace{3mm}
\noindent In other words, the map starts at $\phi(0)=0$, increases by $1/2^N$ on each subinterval which is in $\psi({\cal M}_{c_0,c_1}^N)$, and remains constant on each interval which is not in $\psi({\cal M}_{c_0,c_1}^N)$. Since the sets $\psi({\cal M}_{c_0,c_1}^N)$ are nested, and the corresponding partition $(a_j^N)$ becomes finer with increasing $N$, it follows that the accumulation map becomes lower as $N$ increases. The bottom subplots in Figure~\ref{combs} show, in two separate panels for two different fixed pairs $(c_0,c_1)$, a few instances of  $\phi_{c_0,c_1}^N$, for the same values of $N$ for which the corresponding $\psi({\cal M}_{c_0,c_1}^N)$ sets are plotted in the top subplots.

Consider now the point-wise limit of $\ds \phi_{c_0,c_1} = \lim_{N \to \infty} \phi_{c_0,c_1}^N$. This is also a positive, non-decreasing map of the interval, which we will call the accumulation map of ${\cal M}_{c_0,c_1}$. Some of the topological structure of ${\cal M}_{c_0,c_1}$ for different pairs $(c_0,c_1)$ can be captured by looking at the structure of the corresponding accumulation map $\phi_{c_0,c_1}$, with patterns that may suggest a devil's straircase (which we will discuss below).

Notice that the accumulation maps $\ds \phi_{c_0,c_1}$ and $\ds \phi_{c_0,c_1}^N$ can be expressed as integrals of the indicator function $\chi$ of the subsets ${\cal M}_{c_0,c_1}$ and ${\cal M}_{c_0,c_1}^N$ of the unit interval:
\begin{equation*}
\phi_{c_0,c_1}^N = \int_0^t \chi(\psi({\cal M}_{c_0,c_1}^N))(u) \: du \text{ and respectively } \phi_{c_0,c_1} = \int_0^t \chi(\psi({\cal M}_{c_0,c_1}))(u) \: du.
\end{equation*}

\noindent The structure of the graph of the accumulation map  is specific to the order put on template images in $[0,1]$ (via the map $\psi$). The natural order makes the graph monotonely increasing, so that its maximum value occurs at the end, for $t=1$. This full value of $\phi$ for the pair $(c_0,c_1)$ can be written as the limit $\ds \phi_{c_0,c_1}(1) = \lim_{N \to \infty}\phi_{c_0,c_1}^N(1)$, and coincides with the Lebesgue measure ${\cal L}(\psi({\cal M}_{c_0,c_1}))$. Pulling back, the value of $\phi_{c_0,c_1}(1)$ can be interpreted as the likelihood (with respect to the product measure $\mu $ on $\{0,1 \}^\infty$) for a random template ${\bf s}$ to deliver a bounded critical orbit $o_{\bf s}(0)$ under the iterations of $f_{c_0}$ and $f_{c_1}$ specified by the template ${\bf s}$.  

\begin{defn}
We will say that a pair $(c_0,c_1)$ is $N$-well behaved if the first $N$ iterates of the critical point remain inside the escape disc for \emph{all} $N$-roots, for the given functions $f_{c_0}$ and $f_{c_1}$. We will say that pair is infinitely well behaved if the critical orbit remains bounded under iterations of fixed $f_{c_0}$ and $f_{c_1}$, for almost all templates. 
\end{defn}

\noindent In other words, the pair $(c_0,c_1)$ is $N$-well behaved if the corresponding full $N$-root value $\phi_{c_0,c_1}^N(1)  =1$. The pair is infinitely well behaved if it is $N$-well behaved for all $N$, i.e., if the full template value $\phi_{c_0,c_1}(1) =1$.

\begin{figure}[h!]
\begin{center}
\includegraphics[width=\textwidth]{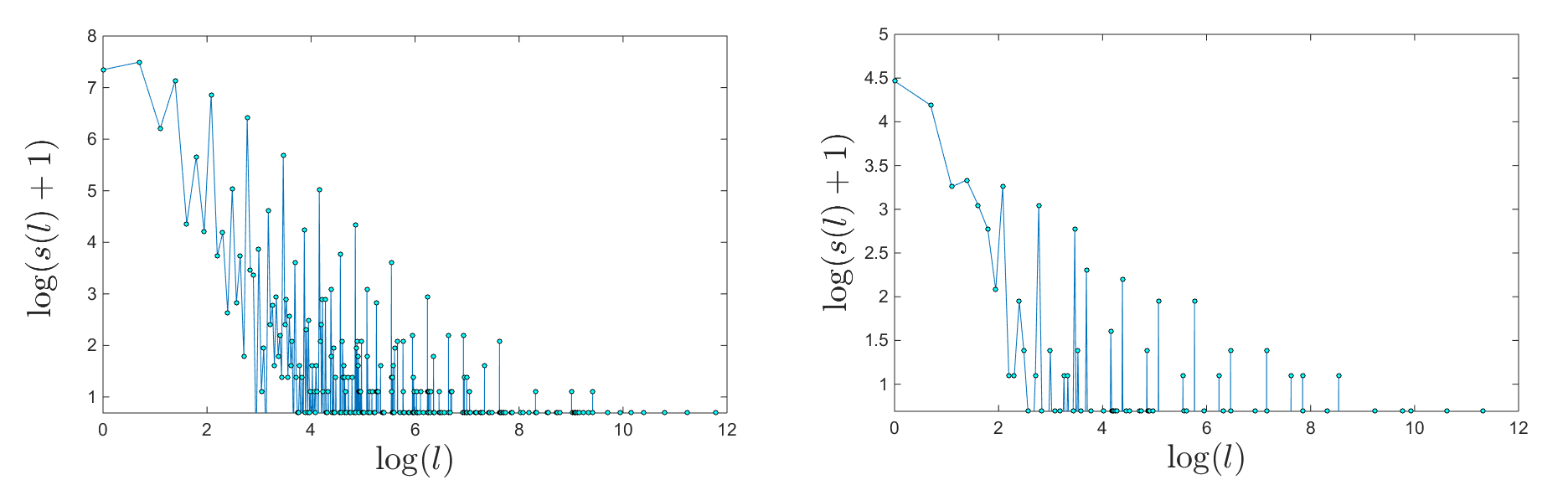}
\end{center}
\caption{\emph{\small {\bf Log-log representation of the distribution of plateau lengths} for $N=20$ and parameters  $c_0=-0.117-0.76i$ and $c_1=-0.62-0.62i$ ({\bf left}); $c_0=-0.75$ and $c_1=-0.117-0.856i$ ({\bf right}).}}
\label{devil}
\end{figure}

To better understand the topological structure of accumulation map graphs, we analyzed the frequency of plateaus of different lengths. In Figure~\ref{devil}, we represent this frequency for the same $(c_0,c_1)$ parameter pairs for which the accumulation maps are shown in Figure~\ref{combs}. For fixed $N=20$, we consider all possible plateau lengths $1/l$, with $1 \leq l \leq 2^N$. For each $l$, we denote by $s(l)$ the number of all plateaux with length $1/l$. A natural question is whether the distribution of plateaux exhibits the power law behavior characteristic to a ``devil's staircase'' (i.e., the number of plateaux of a certain length decreases with the length as a power function, or equivalently as a linear function if the variables are considered in the log-log form). In the two figure panels, we represent $\log(s+1)$ versus $\log(l)$ (one was added to $s$ in order to avoid getting an undefined value for lengths which are not represented). While our log-log plots for finite $N$ appear to be reminiscent of linear behavior (with large variations around a decreasing linear trend), the value ($N=20$) used for these plots is too small to unequivocally extrapolate to linear behavior as $N \to \infty$. In future work, we can explore numerically the structure of these plots for higher values of $N$, as well as for a more complex version of the accumulation map, considering orbits of all critical points (that appear at all steps in the iteration sequence).

A related direction of potential interest is to study the topological structure of the fixed-map Mandelbrot set $\psi({\cal M}_{c_0,c_1})$, as the parameters $c_0$ and $c_1$ are both varied. In the following section, we examine an easier dependence: that of the end value $\phi_{c_0,c_1}(1)$ on the complex parameter pair $(c_0,c_1)$.


\section{Hybrid, contour and multi-Mandelbrot sets}
\label{hybrid_M}

Fix now $c_0 \in \mathbb{D}(2)$. With the notation from Section~\ref{structure}, notice that $b_{c_0}(c_1) = \phi_{c_0,c_1}(1)$, for $c_1 \in \mathbb{C}$. , so that the hybrid Mandelbrot set corresponding to $c_0$ can alternatively be written as:
$${\cal M}_{c_0} = \{ (c_1,\phi_{c_0,c_1}(1)), \text{ for all }  c_1 \in \mathbb{D}(2) \} \subset \mathbb{D}(2) \times [0,1].$$

\noindent One can similarly define the \emph{root hybrid} set corresponding to $c_0$ as
$${\cal M}_{c_0}^N = \{ (c_1,\phi_{c_0,c_1}^N(1)) , \text{ for all }  c_1 \in \mathbb{D}(2) \} \subset \mathbb{D}(2) \times [0,1].$$

\noindent Figure~\ref{hybrid_real}a illustrates the root hybrid set for $c_0=0$, for a root length of $N=20$, with the values of $\phi_{0,c_1}(1)$ corresponding to each $c_1$ represented as colors from blue to dark red. For $c_0=0$, the critical orbit will always be bounded when using the template with all entries zero. With the color map in Figure~\ref{hybrid_real}a, the blue region corresponds to the values of $c_1$ for which $\langle {\bf s} \rangle_{_N} = 0 \hdots 0$ is the only $N$-root that confines the critical orbit to the escape disc. The dark red central region corresponds to the values of $c_1$ for which \emph{all} $N$-roots confine the critical orbit to the escape disc $\mathbb{D}(2)$. In other words, this is the region of $c_1$ for which the pair $(0,c_1)$ is $N$-well behaved. 

Figure~\ref{hybrid_real}b illustrates the central plateau of the hybrid set by comparison with an approximation of the traditional Mandelbrot set ${\cal M}$, computed based on the same number $N$ of iterates. The former is a subset of the latter, a property which remains true in the limit as $N \to \infty$:

\begin{figure}[h!]
\begin{center}
\includegraphics[width=0.8\textwidth]{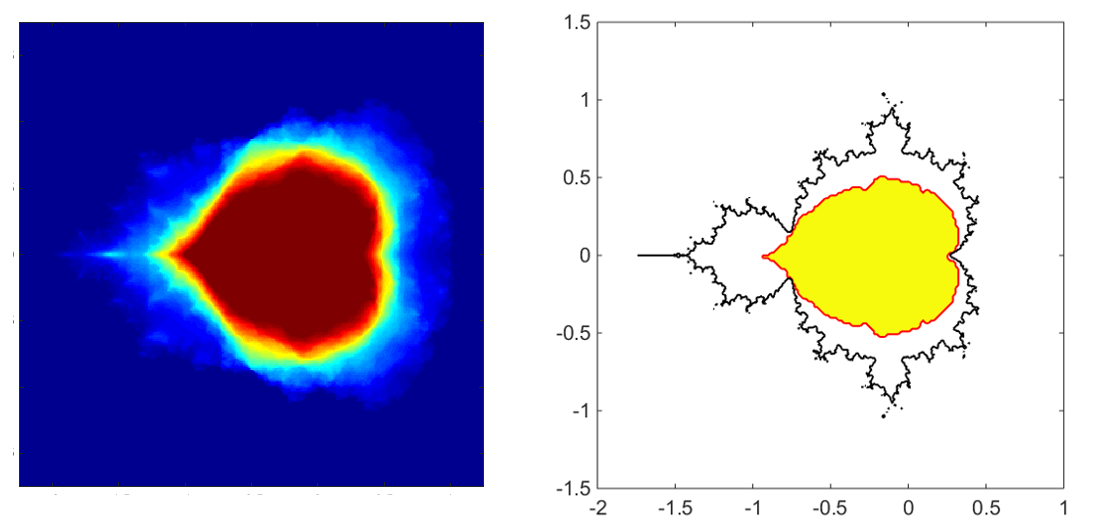}
\end{center}
\caption{\emph{\small {\bf Root hybrid Mandelbrot set ${\cal M}_0^{N}$ for $N=20$.} {\bf Left.} Each pair $(c_1,b)$ is represented as a point $c_1$  in the complex plane with an associated color from blue to dark red for $b=\phi_{0,c_1}^{N}(1)$ increasing from 0 to 1. {\bf Right.} The shaded region corresponds to the $N$-well behaved pairs $\phi_{0,c_1}^N(1)=1$ for $N=20$. It is illustrated as a subset of a truncated representation of the traditional Mandelbrot set, based for comparison on $N=20$ iterations as well (the interior of the black curve).}}
\label{hybrid_real}
\end{figure}

$$\{ c_1 \in \mathbb{D}(2) \text{ with } \phi_{0,c_1}(1)=1 \} \subset {\cal M} \subset \mathbb{D}(2).$$

\begin{thm}
The set of infinitely well-behaved pairs $(0,c_1)$ contains the disc centered at the origin, of radius $1/4$.
\end{thm}

\noindent {\bf Proof.} Consider $c_1 \in \mathbb{C}$, with $\lvert c_1 \rvert < 1/4$. Then $\Delta = 1-4 \lvert c_1 \rvert >0$, and one can define $\ds d = \frac{1 +\sqrt{\Delta}}{2} = \frac{1+\sqrt{1-4\lvert c_1 \rvert}}{2}$. Notice that $0<d<1$. Suppose now that $\lvert z \rvert <d$. It follows that 
\begin{eqnarray*}
\lvert z^2 + c_1 \rvert &\leq& \lvert z \rvert^2 + \lvert c_1 \rvert < d^2 + \lvert c_1 \rvert = \frac{1+2\sqrt{\Delta} + \Delta}{4} + \lvert c_1 \rvert = \frac{1+\sqrt{\Delta}}{2} = d.
\end{eqnarray*}

\noindent Subsequently, if $\lvert z \rvert <d$ then $\lvert f_{c_{{\bf s}_n}}(z) \rvert<d$, for all $n \geq 1$. Inductively, it follows that the orbit of $\xi_0=0 < d$ is completely contained in the disc of radius $d$ around the origin, hence it is bounded. In conclusion: the critical orbit is bounded for all possible template iterations under the pair of maps $(0,c_1)$ with $\lvert c_1 \rvert < 1/4$. 

\hfill $\Box$ 

%

\noindent Based on our simulations, we conjecture the following:

\begin{conj}
The set of all $c_1 \in \mathbb{D}(2)$ such that the pair $(0, c_1)$ is infinitely well-behaved, is a connected subset of ${\cal M}$.
\end{conj}

\begin{figure}
\begin{center}
\includegraphics[width=0.6\textwidth]{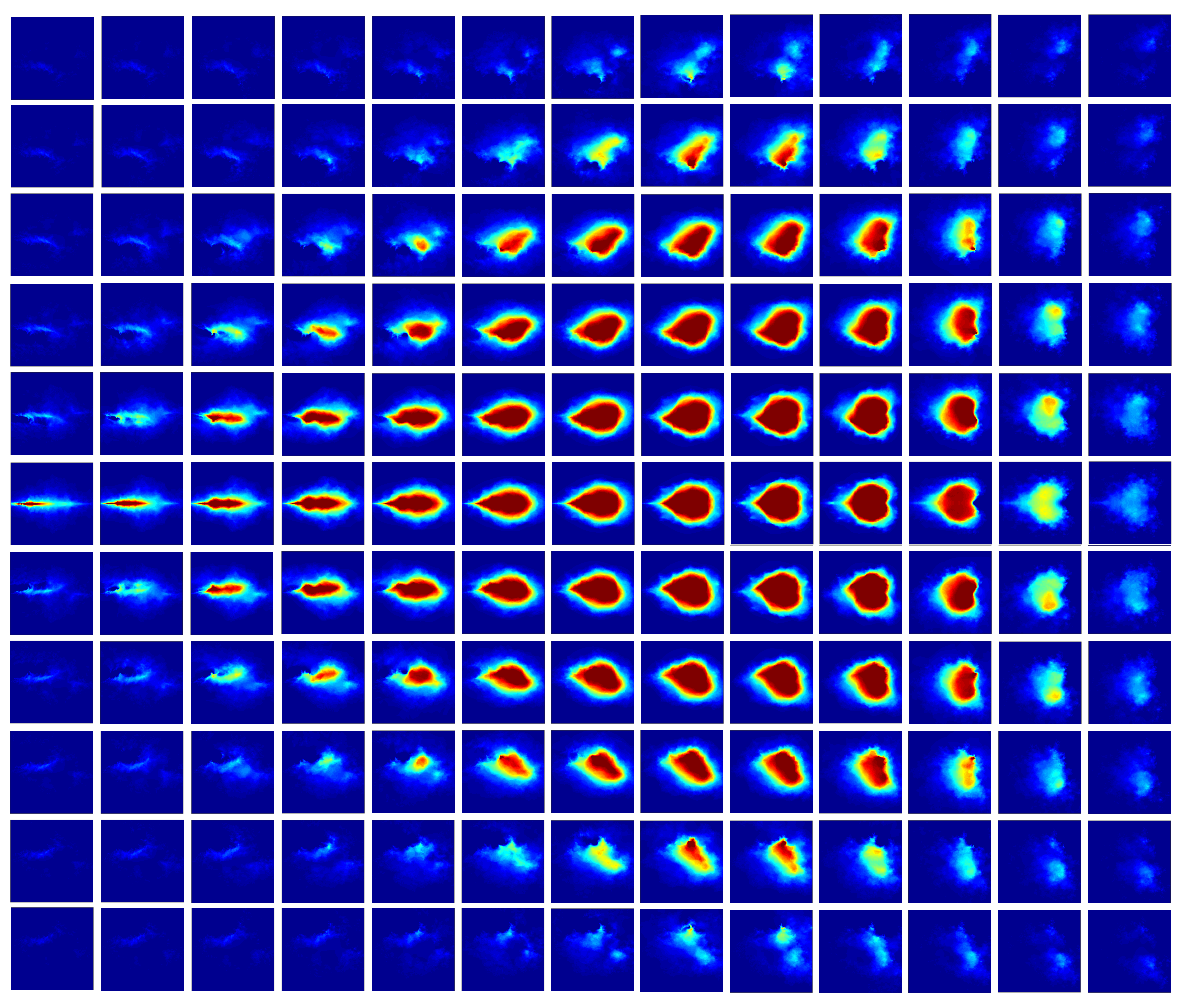}
\end{center}
\caption{\emph{\small {\bf Hybrid sets for a grid of $c_0$ values.} Each panel represents one value of $c_0$, covering the interval $[-1.6, 0.8]$ along the real axis and $[-1, 1]$ along the imaginary axis, with distance $0.2$ in between sample values in both directions.}}
\label{hybrid_table}
\end{figure}

\noindent While the hybrid set for $c_0=0$ may be the easiest to study and understand, one can consider hybrid sets for other values of $c_0$. Figure~\ref{hybrid_table} illustrates the hybrid sets for a grid of $c_0$ complex values with Re$(c_0) \in [-1.6,0.8]$, and Im$(c_0) \in [-1,1]$. There are a few different ways in which one can organize this atlas of $c_0$-indexed hybrid sets. Below we describe two such ways, by defining two new sets. Based upon the function $b$ defined in Section~\ref{intro}, one can define a new function 
\begin{eqnarray*}
\beta \colon \mathbb{D}(2) \to [0,1] \text{, given by } \beta(c_0) &=& \text{max}\{b(c_0,c_1), \text{ for } c_1 \in \mathbb{D}(2) \}\\ 
&=& \text{max}\{\phi_{c_0,c_1}(1), \text{ for } c_1 \in \mathbb{D}(2) \}.
\end{eqnarray*}

\noindent We also define a variation for truncated templates,
$$\beta^N \colon \mathbb{D}(2) \to [0,1] \text{, given by } \beta^N(c_0) = \text{max}\{\phi^N_{c_0,c_1}(1), \text{ for } c_1 \in \mathbb{D}(2) \}.$$

\begin{defn}
The contour Mandelbrot set is the graph of $\beta$:
$${\cal CM} = \{ (c_0,\beta(c_0)), \text{ for all } c_0 \in \mathbb{D}(2) \} \subset \mathbb{D}(2) \times [0,1].$$
Similarly, the $N$-root contour Mandelbrot set is
$${\cal CM}^N = \{ (c_0,\beta^N(c_0)), \text{ for all } c_0 \in \mathbb{D}(2)  \} \subset \mathbb{D}(2) \times [0,1].$$
\end{defn}

\begin{figure}[h!]
\begin{center}
\includegraphics[width=0.8\textwidth]{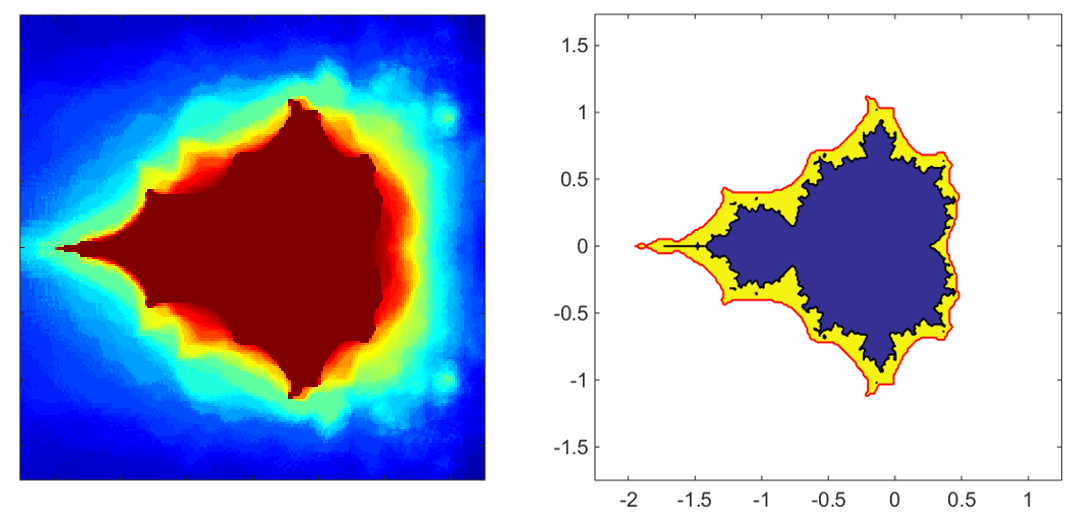}
\end{center}
\caption{\emph{\small {\bf $N$-root contour Mandelbrot set,} for $N=8$. {\bf Left.} The colors represent values of $b$ from zero (blue) to 1 (dark red). {\bf Right.} The yellow region represents the central plateau of ${\cal CM}^N$ for $N=20$. In blue is shown an approximation of the classical Mandelbrot set ${\cal M}$, for the same $N=20$.}}
\label{mega_hybrid}
\end{figure}

\noindent In other words, the contour Mandelbrot set assigns to every $c_0 \in \mathbb{D}(2)$ the largest value of $\phi_{c_0,c_1}(1)$ over all $c_1 \in \mathbb{D}(2)$. Similarly, one can assign to every $c_0$ the largest value of $\phi_{c_0,c_1}^N(1)$, and obtain the $N$-root contour Mandelbrot set -- illustrated in Figure~\ref{mega_hybrid}a for $N=8$, with the colors representing the level sets of the function $\beta$.

One can consider the level set ${\cal PM}^N$ corresponding to the highest value of $\beta$, illustrated in Figure~\ref{mega_hybrid}a as the central (burgundy) plateau of ${\cal CM}^N$
\begin{equation*}
{\cal PM}^N = \{ c_0 \in \mathbb{D}(2) \text{ for which there exists a } c_1 \in \mathbb{D}(2) \text{ such that } (c_0,c_1) \text{ is  $N$-well behaved} \}.
\end{equation*}
 
 \noindent One can also consider the $N$-truncated Mandlebrot set
\begin{equation*}
{\cal M}^N = \{ c_0 \in \mathbb{D}(2) \text{ such that } f_{c_0}^{\circ k}(0) \in \mathbb{D}(2) \text{ for all } 0 \leq k \leq N \}.
\end{equation*} 
  
\noindent It is easy to see that, if $c_0 \in {\cal M}^N$, then there exists $c_1=c_0$, for which $(c_0,c_0)$ is $N$ well behaved, that is $c_0 \in {\cal PM}^N$. Hence ${\cal M}^N \subset {\cal PM}^N$, as illustrated in Figure~\ref{mega_hybrid}b. It follows, in the $N \to \infty$ limit, that 
\begin{equation*}
{\cal M} \subset \{ c_0 \in \mathbb{D}(2) \text{ for which there is a }c_1 \in \mathbb{D}(2) \text{ such that } (c_0,c_1) \text{ is infinitely well behaved} \}.
\end{equation*}

\begin{defn}
The multi-Mandelbrot set is defined as
$${\cal MM} = \{ (c_0,c_1) \in \mathbb{D}^2(2), \text{ such that } \phi_{c_0,c_1}(1) =1 \}.$$

\noindent and the $N$-root multi-Mandelbort set is
$${\cal MM}^N = \{ (c_0,c_1) \in \mathbb{D}^2(2), \text{ such that } \phi_{c_0,c_1}^N(1) =1 \}.$$
\end{defn}

\begin{figure}[h!]
\begin{center}
\includegraphics[width=0.9\textwidth]{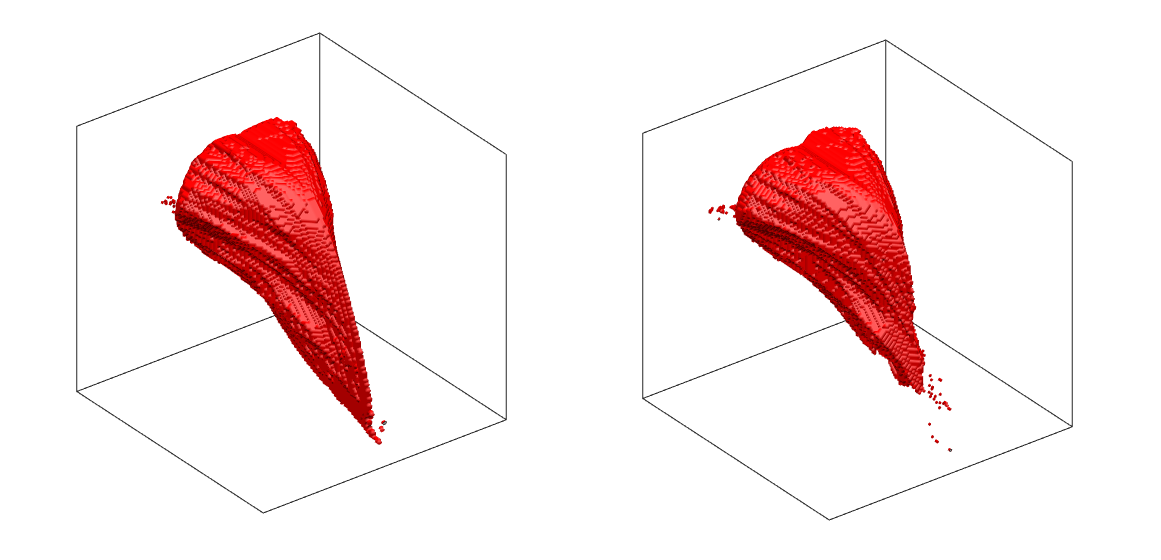}
\end{center}
\caption{\emph{\small {\bf Three-dimensional slices of the $N$-root multi-Mandelbrot set,} for $N$=8. {\bf Left.} Real $c_0 \in [-2,1]$. {\bf Right.} Complex $c_0= \text{Re}(c_0) + 0.1i$, with $\text{Re}(c_0) \in [-2,1]$.}}
\label{super_hybrid}
\end{figure}

\noindent ${\cal MM}$ is the set of all infinitely well behaved parameter pairs in $\mathbb{D}(2) \times \mathbb{D}(2)$. Both ${\cal MM}$ and ${\cal MM}^N$ have real dimension 4. In Figure~\ref{super_hybrid}, we illustrate 3-dimensional slices of ${\cal MM}^N$. While Figure~\ref{super_hybrid} shows that not all 3D slices of ${\cal MM}^N$ are connected, we conjecture that:

\begin{conj}
${\cal MM}$ is a connected set in $\mathbb{C}^2$.
\end{conj}

\section{Multicritical definitions}
\label{multicrit}

\noindent  Our condition for the initial critical point at zero to have a bounded orbit, while inspired by its relevance to potential applications, differs from the postcritically bounded condition typically used in random iterations~\cite{comerford2013preservation}. In order to reconcile our results with this framework, we investigate whether the same objects can be studied when using the alternative definition, in which \emph{all} critical points, generated at all stages of the iteration, are required to remain bounded (we will call the corresponding sets {\bf multi-critical M-sets}).

This idea follows naturally when one considers the significance of a critical point in the broader context of non-autonomous iterations. When one iterates consecutively $k$ differentiable functions $f_1$ to $f_k$, having a critical point for the composition requires that
$$(f_k \circ f_{k-1} \circ \hdots \circ f_2 \circ f_1)'(z) = f_k'[f_{k-1} \circ \hdots \circ f_2 \circ f_1(z)] \cdot \hdots  \cdot f_2'[f_1( z)] \cdot f_1'(z) =0.$$

\noindent Hence $z$ is a critical point for the composition if and only if either $z$ is a critical point for $f_1$ or $f_{j-1} \circ \hdots \circ f_2 \circ f_1(z)$ is a critical point for $f_j$, for some $2 \leq j \leq k$. In the case of template iterations, we are applying a sequence of quadratic functions $f_{c_0}(z)$ and $f_{c_1}(z)$. The situation is quite simple, because, at each step, the function we are iterating has only one critical point $\xi_0=0$ (independently on the step). So an efficient computational strategy for determining whether a template iteration is postcritically bounded consists of checking whether the tails of the orbits below are bounded or whether they escape to $\infty$:\\

$0 \stackrel{f_{c_{s_1}}}{\longrightarrow} * \stackrel{f_{c_{s_2}}}{\longrightarrow} * \stackrel{f_{c_{s_3}}}{\longrightarrow} * \stackrel{f_{c_{s_4}}}{\longrightarrow} * \hdots \stackrel{f_{c_{s_k}}}{\longrightarrow} * \stackrel{f_{c_{s_{k+1}}}}{\longrightarrow} \hdots$ 

$\hspace{1cm} 0 \stackrel{f_{c_{s_2}}}{\longrightarrow} * \stackrel{f_{c_{s_3}}}{\longrightarrow} * \stackrel{f_{c_{s_4}}}{\longrightarrow} * \hdots \stackrel{f_{c_{s_k}}}{\longrightarrow} * \stackrel{f_{c_{s_{k+1}}}}{\longrightarrow} \hdots$ 

$\hspace{2cm} 0 \stackrel{f_{c_{s_3}}}{\longrightarrow} * \stackrel{f_{c_{s_4}}}{\longrightarrow} * \hdots \stackrel{f_{c_{s_k}}}{\longrightarrow} * \stackrel{f_{c_{s_{k+1}}}}{\longrightarrow} \hdots$ 

\hspace{5cm}$\vdots$

$\hspace{3cm} 0  \stackrel{f_{c_{s_4}}}{\longrightarrow} * \hdots \stackrel{f_{c_{s_k}}}{\longrightarrow} * \stackrel{f_{c_{s_{k+1}}}}{\longrightarrow} \hdots$

$\hspace{4.5cm} 0 \stackrel{f_{c_{s_k}}}{\longrightarrow} * \stackrel{f_{c_{s_{k+1}}}}{\longrightarrow} \hdots$\\

\noindent Hence we consider the following alternative definitions, which consider all critical points generated throughout the iteration:

\begin{defn}
Fix a template ${\bf s} \in \{ 0,1 \}^\infty$. The corresponding \textbf{multi-critical fixed-template Mandelbrot set} is defined as
$$^m{\cal M}_{\bf s} = \{ (c_0,c_1) \in \mathbb{D}(2)^2, \text{ where } o_{{\bf s}^k}(0) \text{ is bounded for all integers } k \geq 0  \}.$$
\noindent where ${\bf s}^k$ is the right $k$-shift of the template ${\bf s}$. We additionally define the $N$-root multi-critical fixed-template Mandelbrot set as
$$^m{\cal M}_{\bf s}^N = \{ (c_0,c_1) \in \mathbb{D}(2)^2, \text{ where }  \lvert o_{{\bf s}^k}^{N-k}(0) \rvert \leq 2 \text{ for all } 1 \leq k \leq N  \}.$$
\end{defn}

\vspace{2mm}
\begin{defn}
Fix $(c_0,c_1) \in \mathbb{D}(2)^2$. The corresponding \textbf{multi-critical fixed-map Mandelbrot set} is defined as
$$^m{\cal M}_{c_0,c_1} = \{ {\bf s} \in \{ 0,1 \}^\infty, \text{ where } o_{{\bf s}^k}(0) \text{ is bounded for all integers } k \geq 0 \}.$$

\noindent The $N$-root multi-critical fixed-map Mandelbrot set is defined as
$$^m{\cal M}^N_{c_0,c_1} = \{ {\bf s} \in \{ 0,1 \}^\infty, \text{ where }  \lvert o_{{\bf s}^k}^{N-k}(0) \rvert \leq 2 \text{ for all } 1 \leq k \leq N \}.$$
\end{defn}

\noindent We consider the multi-critical version of the function $b$ defined in Section~\ref{intro}, as well as its variation for truncated templates:
$$^m b \colon \mathbb{C}^2 \to [0,1], \: \text{ given by } ^m b(c_0,c_1) = {\cal L} \left( \psi(^m {\cal M}_{c_0,c_1}) \right),$$
$$^m b^N \colon \mathbb{C}^2 \to [0,1], \: \text{ given by } ^m b^N(c_0,c_1) = {\cal L} \left( \psi(^m {\cal M}^N_{c_0,c_1}) \right).$$

\noindent Then, for a fixed $c_0 \in \mathbb{C}$, we define the projections
$$^m b_{c_0} \colon \mathbb{C} \to [0,1], \: \text{ given by } ^m b_{c_0}(c_1) = ^m b(c_0,c_1),$$
$$^m b^N_{c_0} \colon \mathbb{C} \to [0,1], \: \text{ given by } ^m b_{c_0}(c_1) = ^m b(c_0,c_1).$$

\vspace{2mm}
\begin{defn}
For a fixed $c_0 \in \mathbb{D}(2)$, the \textbf{multi-critical hybrid Mandelbrot set} is defined as:
$$^m{\cal M}_{c_0} = \{ (c_1,^m b_{c_0}(c_1)) \in \mathbb{D}(2) \times [0,1], \text{ for all } c_1 \in \mathbb{D}(2) \}$$
\noindent We additionally define the $N$-root multi-critical hybrid Mandelbrot set as:
$$^m{\cal M}_{c_0}^N = \{ (c_1,^m b^N_{c_0}(c_1)) \in \mathbb{D}(2) \times [0,1], \text{ for all } c_1 \in \mathbb{D}(2) \}$$
\end{defn}

\begin{figure}[h!]
\begin{center}
\includegraphics[width=0.9\textwidth]{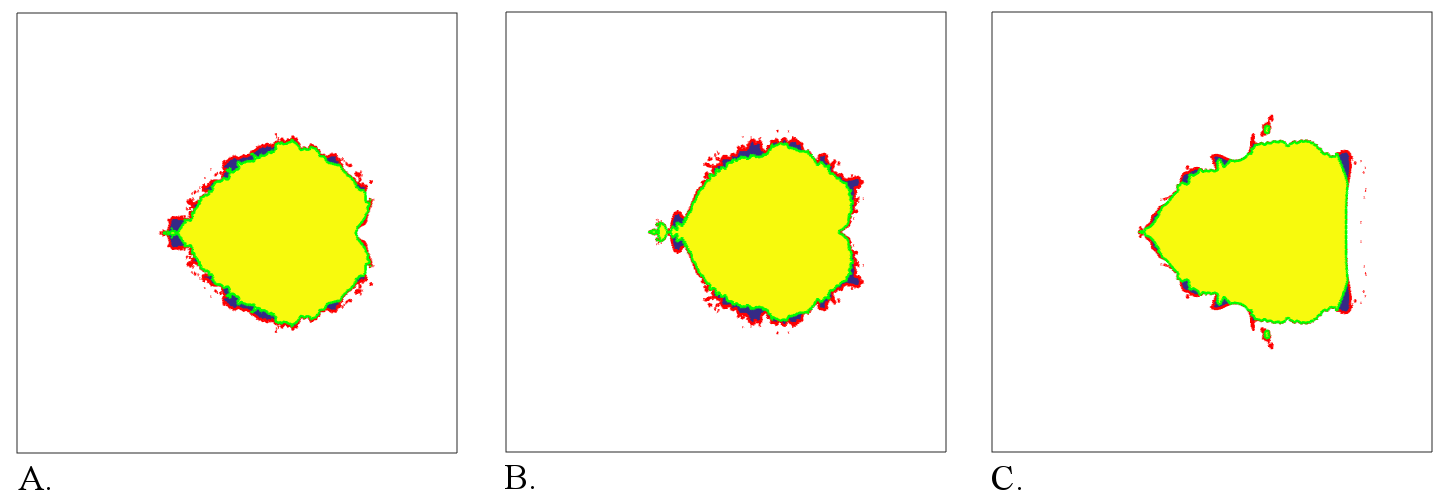}
\end{center}
\caption{\emph{\small {\bf Examples of fixed-template M-slices and multi-critical M-slices for $c_0=0$.} The template ${\bf s}$ was selected at random such that {\bf A.} the probability of a 1 entry is $p=0.5$; {\bf B.} the probability of a 1 entry is $p=0.25$; {\bf C.} the probability of a 1 entry is $p=0.75$. The larger region (red contour, shaded in blue) represents the M-slice for $c_0=0$, which is the subset $c_1 \in \mathbb{C}$ for which the original critical point $\xi_0=0$ is bounded under the maps $(c_0,c_1)$ with template ${\bf s}$.  The nested region (green contour, shaded yellow) represents the multi-critical M-slice, which is the subset of $c_1 \in \mathbb{C}$ for which the whole critical set (consisting of $\xi_0=0$ from all iteration steps, as described in the text) is bounded under the iteration system.}}
\label{Mset_multicrit}
\end{figure}

\vspace{2mm}
\noindent  Figure~\ref{Mset_multicrit} shows a few examples of multi-critical M-slices for fixed $c_0=0$ and for fixed templates with specific frequencies of iterating one function versus the other. These sets are presented in comparison with their regular M-slice counterparts, of which they are subsets. 

\begin{figure}[h!]
\begin{center}
\includegraphics[width=0.8\textwidth]{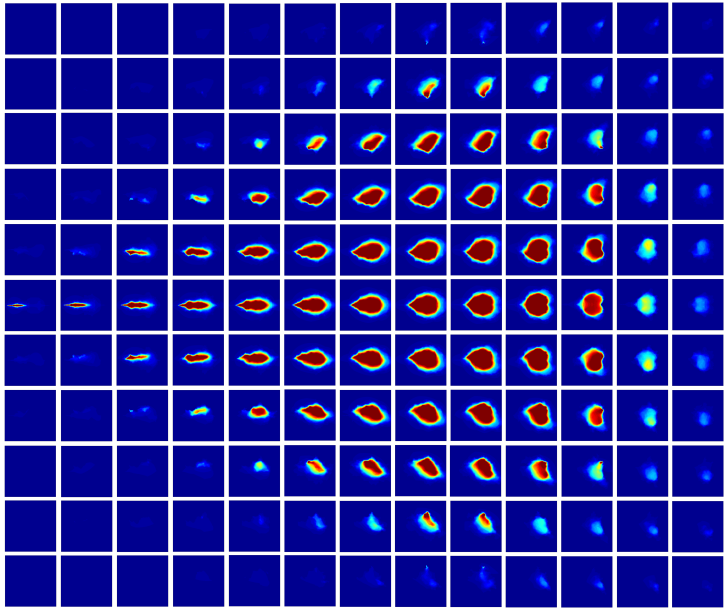}
\end{center}
\caption{\emph{\small {\bf Hybrid sets for a grid of $c_0$ values}. Regular hybrid sets (left) and multi-critical hybrid sets (right). The grid covers $[-1.6, 0.8]$ along the real axis and $[-1, 1]$ along the imaginary axis, with distance $0.2$ in between sample values in both directions. All panels are based on $N$-root approximations with $N=8$ iterations.}}
\label{grids}
\end{figure}

As in the previous section, we obtain a more comprehensive view by computing multi-critical \emph{hybrid} sets. For illustration purposes, we consider the same parameter values $c_1 \in \mathbb{C}$ as in Figure~\ref{hybrid_table} in Section~\ref{hybrid_M}, and show in Figure~\ref{grids}, side by side, the tables for regular hybrid slices and for multi-critical hybrid slices. Notice that the burgundy plateaux are identical between the two figures, for all $c_0$ panels (as justified in the following section). Moreover, the regular and multi-critical hybrid sets differ less when $c_0$ is close to the origin, and more in the panels for $c_0$ away from the origin.

\subsection{Fatou-Julia theorem}
\label{connectedness}

A great rationale for using the multi-critical definition for the Mandelbrot sets is the relationship of this condition with the topology of the template Julia set. IT is a well-known phenomenon in both real and complex dynamics that the behavior of a polynomial's critical set encompasses the whole dynamic behavior of the map. In particular within the complex quadratic family $f_c(z)=z^2+c$ (with a unique critical point at $z_0=0$ for all maps), a bounded orbit for $z_0=0$ (i.e., $c$ in the Mandelbrot set) is equivalent to connectedness of the Julia set $J(f_c)$, and an escaping orbit for $z_0=0$ is equivalent to the Julia set being totally disconnected. This duality stands more generally for polynomials $P$ of degree $d \geq 2$, where the filled Julia set is connected (a cellular set, in fact) if and only if it contains all finite critical points of $P$, and it has uncountably many connected components iff the orbit of at least one of the critical points escapes (see, for example Theorem 9.5 in~\cite{milnor2006dynamics}, or Theorem 4.1 in~\cite{carleson1993complex}).

\begin{figure}[h!]
\begin{center}
\includegraphics[width=0.6\textwidth]{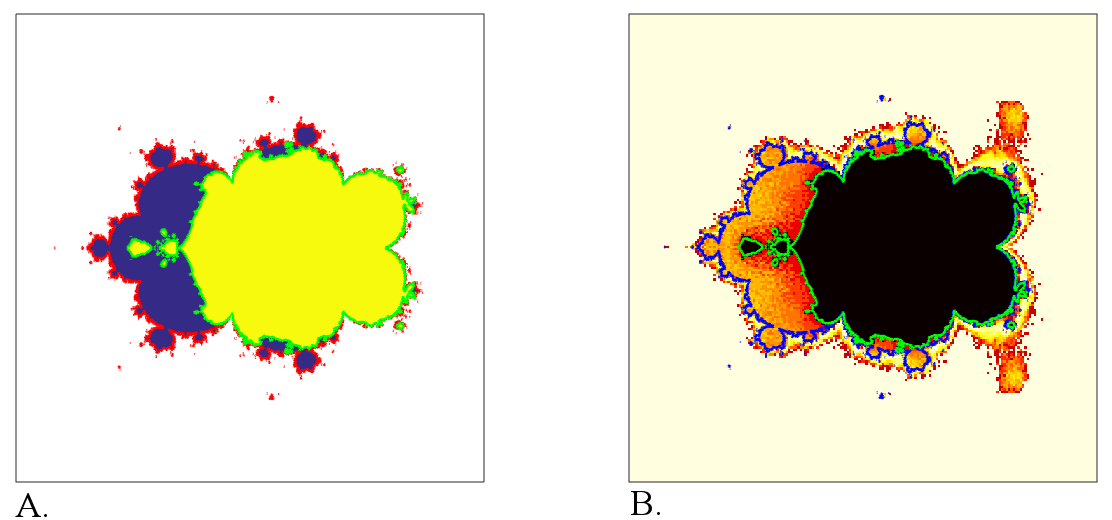}
\caption{\emph{\small {\bf Mandelbrot slice, multicritical Mandelbrot slice and connectedness slice for fixed periodic template ${\bf s} = [011]$ and fixed $c_0=0$.} {\bf A.} The Mandelbrot slice is the blue interior of the red curve. The multicritical Mandelbrot slice is the yellow interior of the green curve, which is a subset of the former. {\bf B.} The panel illustrates the connectedness locus of the Julia set in the slice $c_0=0$, with colors in the hot spectrum representing the number of connected components of the Julia set, estimated numerically. Black represents the parameters $c_1$ for which the template Julia was estimated to have a single connected component (the connectedness locus). The intermediate darker to lighter (brown to yellow) colors represent the increasing number of connected components (estimated by our algorithm to be finite, but $\geq 2$). The outside (white) region represents the locus where our algorithm estimated the Julia set to be totally disconnected. The boundary of the corresponding Mandelbrot and muticritical Mandelbrot slices (shown in {\bf A}) are overlaid as blue and respectively green curves, for comparison.}}
\label{connected_periodic}
\end{center}
\end{figure}

We want to check how this duality holds in the context of template iterations. Clearly, the result follows trivially in the case of periodic templates, since the system is equivalent in this case with the iteration of a single polynomial of higher degree, with critical orbits identical with those initiated at all steps of the template iteration. We illustrate this in Figure~\ref{connected_periodic}, for the periodic template ${\bf s} =[011]$, showing that the multi-crititical Mandelbrot set coincides in this case with the connectedness locus of the Julia set in $\mathbb{C}$.

\begin{figure}[h!]
\begin{center}
\includegraphics[width=0.6\textwidth]{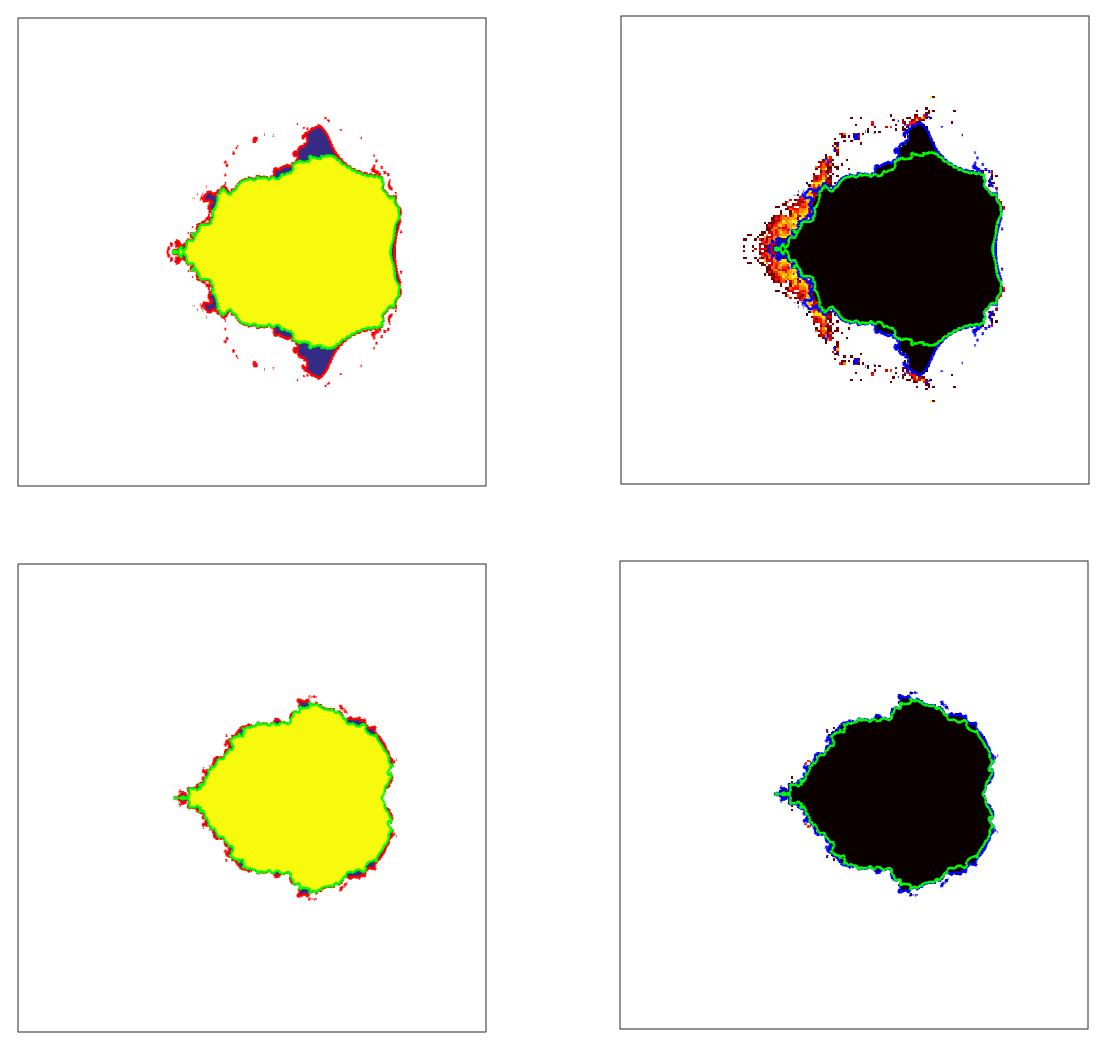}
\caption{\emph{\small {\bf Mandelbrot slice, multicritical Mandelbrot slice and connectedness slice for $c_0=0$ and a random template}. The top panels illustrate these objects using an approximation based on a 50 entry root of the template; the bottom panels show a refinement, based on a longer root of 200 entries of the same template. {\bf Left.} In each case, the Mandelbrot slice is the blue interior of the red curve. The multicritical Mandelbrot slice is the yellow interior of the green curve, which is a subset of the former. {\bf Right.} The panels illustrate the connectedness of the Julia set in the slice $c_0=0$, with colors in the hot (red to yellow) spectrum representing the number of connected components of the Julia set, estimated numerically. Black represents the parameters associated with one connected component (the connectedness locus). As before, the intermediate darker to lighter (brown to yellow) colors represent the increasing number of connected components (estimated by our algorithm to be finite, but $\geq 2$). The outside (white) region represents the locus where our algorithm estimated the Julia set to be totally disconnected. The intermediate lighter to darker (yellow to brown) colors represent the number of connected components estimated by our algorithm to be finite, but $\geq 2$. The boundary of the corresponding Mandelbrot and muticritical Mandelbrot slices (shown on the left) are overlaid as blue and respectively green curves, for comparison.}}
\label{connected_random}
\end{center}
\end{figure}

For random iterations, the traditional dichotomy of the Julia sets having either one or uncountably many connected components is no longer guaranteed, since Julia sets may exist that have a finite number $\geq 2$ of connected components. However, the equivalence between the postcritically bounded locus and the Julia set connectedness locus remains valid, even for templates which are not periodic (as illustrated in Figure~\ref{connected_random}). That is because, for any $m>0$, the Green's function defined by ~\cite{fornaess1991random}
$$G_m(z) = \lim_{n \to \infty} \frac{1}{2^n} \log \lvert Q_{m,m+n+1}(z) \rvert$$
\noindent still exists on the iterated basin of infinity ${\cal A}_{\infty,m} = \{ z \in \mathbb{C} \text{ such that } \lim \lvert Q_{m,n}(z) \rvert =\infty \}$ (where $Q_{m,n}(z) = f_{{\bf s}_{n}} \circ \hdots \circ f_{{\bf s}_{m+1}}$), so that $o_{\bf s}(z)$ escapes iff $z \in {\cal A}_{\infty,m}$ iff $G_m(z) >0$. In addition, it was also shown that there is a B\"{o}ttcher isomorphism $\varphi_m \colon {\cal A}_{\infty,m} \to \overline{\mathbb{C}} \setminus \overline{\mathbb{D}}$ such that $\varphi_{m+1} \circ f_{{\bf s}_{m+1}} = \varphi_m(z^2)$~\cite{bruck2001geometric}. The result follows similarly to the traditional result in the case of single map iterations.

\subsection{Hybrid sets}

Section~\ref{connectedness} shows that, at the level of each specific template iteration (for fixed parameters $c_0$ and $c_1$, and a fixed template ${\bf s}$), the multi-critical Mandelbrot set contains more information about the system than the regular Mandelbrot set. However, for practical purposes, a regular Mandelbrot set is computationally less expensive. In mind with the goal of obtaining ``averaged'' information over all templates (useful for applications), we want to investigate whether one definition remains superior to the other in the context of hybrid Mandelbrot sets.

As before, we can define the functions
$$^m \beta \colon \mathbb{D}(2) \to [0,1] \text{, given by } ^m \beta(c_0) = \text{max}\{^m b(c_0,c_1), \text{ for } c_1 \in \mathbb{D}(2) \} \text{ and }$$
$$^m \beta^N \colon \mathbb{D}(2) \to [0,1] \text{, given by } ^m \beta^N(c_0) = \text{max}\{^m b^N(c_0,c_1), \text{ for } c_1 \in \mathbb{D}(2) \}.$$

\vspace{3mm}
\begin{defn}
The multicritical contour Mandelbrot set is the graph of $^m \beta$:
$${\cal CM} = \{ (c_0,^m \beta(c_0)), \text{ for all } c_0 \in \mathbb{D}(2) \} \subset \mathbb{D}(2) \times [0,1].$$
Similarly, one can define the $N$-root multicritical contour Mandelbrot set as
$${\cal CM}^N = \{ (c_0,^m \beta^N(c_0)), \text{ for all } c_0 \in \mathbb{D}(2)  \} \subset \mathbb{D}(2) \times [0,1].$$
\end{defn}

\begin{figure}[h!]
\begin{center}
\includegraphics[width=0.7\textwidth]{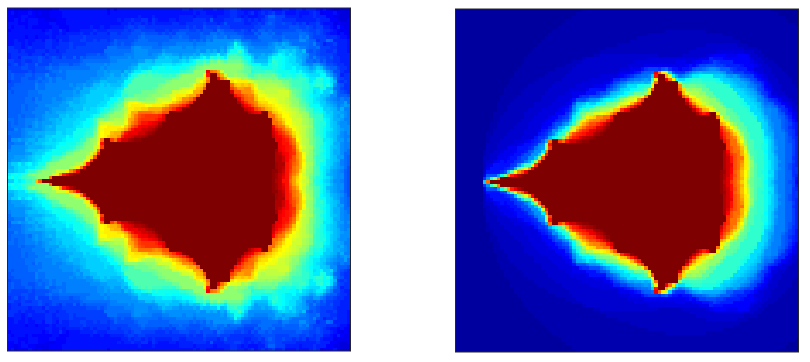}
\end{center}
\caption{\emph{\small {\bf Contour Mandelbrot sets} for the regular (left) and multi-critical (right) definitions. For each point $c_0 \in \mathbb{C}$, the color represents the value of the maximum likelihood $b$, in the space of templates of length $N=8$, for the critical orbit(s) to be bounded, with the maximum taken over all $c_1 \in \mathbb{C}$. The color map ranges from zero (blue) to 1 (dark red). Improvements to the quality of the illustrations can be obtained with increasing the spatial resolution in the complex $c_0$ plane, and the number $N$ of iterations, all of which increase multiplicatively the computational cost.}}
\label{centrals}
\end{figure}

\vspace{2mm}
\begin{defn}
We define the multi-critical multi-Mandelbrot set as
$$^m{\cal MM} = \{ (c_0,c_1) \in \mathbb{D}(2)^2, \text{ such that } o_{{\bf s}^k}(0) \text{ is bounded for almost all templates } {\bf s} \text{ and all } k \geq 0 \}$$
\noindent where ${\bf s}^k$ designates the $k$-shift of ${\bf s}$. We define the $N$-root multi-critical multi-Mandelbrot set as
$$^m{\cal MM}^N =  \{ (c_0,c_1) \in \mathbb{D}(2)^2, \text{ with } \lvert o_{{\bf s}^k}^{N-k}(0) \rvert \leq 2 \text{ for almost every template } {\bf s} \text{ and all } 0 \leq k \leq N \}$$
\end{defn}

\vspace{2mm}
\noindent Recall that:
$${\cal MM}^N = \{ (c_0,c_1) \in \mathbb{D}^2(2), \text{ such that } \phi_{c_0,c_1}^N(1) =1 \}$$

Suppose $(c_0,c_1) \in \, ^m{\cal MM}^N$, that is $\lvert o_{\langle {\bf s}^k \rangle_{N-k}}(0) \rvert \leq 2$ for almost every ${\bf s} \in \{ 0,1 \}^\infty$ and all shifts $0 \leq k \leq N$. This is true in particular for $k=0$, which implies that $\lvert o_{\langle {\bf s} \rangle_N}(0) \rvert \leq 2$ for almost every ${\bf s}$, hence $(c_0,c_1) \in \, {\cal MM}^N$. This proves that $^m{\cal MM}^N \subseteq {\cal MM}^N$.

Now suppose $(c_0,c_1) \in \, {\cal MM}^N$, This means that $o_{\langle {\bf s} \rangle_N}(0)$ is bounded under $(c_0,c_1)$ iterations for almost every ${\bf s} \in \{ 0,1 \}^\infty$, hence $\lvert o_{\langle {\bf s} \rangle_N}(0) \rvert \leq 2$ for every ${\bf s} \in \{ 0,1 \}^\infty$. The inverse image under the shift map $\sigma$ of a subset of full measure in $\{ 0, 1 \}^\infty$  has full measure as well. In addition, a finite intersection of sets of full measure has itself full measure. Hence ${\cal MM}^N \subseteq \, ^m{\cal MM}^N$, and thus $^m{\cal MM}^N = {\cal MM}^N$ for all $N$. Since
$${\cal MM}= \bigcap_{n=1}^\infty {\cal MM}^N \text{ and }  ^m{\cal MM} =  \bigcap_{n=1}^\infty \, ^m{\cal MM}^N$$

\noindent it follows easily that:

\begin{prop}
For any $(c_0,c_1) \in \mathbb{D}(2)^2$, the multi-Mandelbrot and multi-critical multi-Mandelbrot sets coincide, i.e. ${\cal MM}=$ $^m{\cal MM}$.
\end{prop}

\noindent It also follows as a consequence that, while the contour Mandelbrot sets differ between the two (regular and multi-critical) definitions, their central plateaux are identical (see Figure~\ref{centrals}):

\begin{prop}
For any $(c_0,c_1) \in \mathbb{D}(2)^2$, the central plateaux of ${\cal CM}$ and $^m{\cal CM}$ coincide.
\end{prop}

\noindent It is useful that the two definitions are interchangeable at the level of these averaged sets. One of the advantages is the possibility to compute these sets using the more efficient definition (with one critical point), yet obtain the connectedness of the Julia set delivered by the stricter (multi-critical) definition, as discussed in Section~\ref{connectedness}. In Section~\ref{discussion}, we further discuss the relationships between the two means of computing these sets, and we interpret this in the context of applications.

\section{Discussion}
\label{discussion}

\subsection{Interpretation of our results}

In this paper, we expanded our study of asymptotic dynamics under template iterations, introducing new concepts and ideas for further research. Inspired by one of the alternative definitions of the traditional Mandelbrot set for single map iterations, we investigated extensions in the new context of template iterations. Since in this case the parameters are a combination of the complex value pair $(c_0,c_1) \in \mathbb{C}^2$ and the template ${\bf s} \in \{ 0,1 \}^\infty$, we defined a few different types of parameter slices, which can be more easily visualized and understood than the whole parameter set $\mathbb{C}^2 \times \{ 0,1 \}^\infty$.

We first embraced an application-driven idea, and viewed the Mandelbrot set as the parameter locus for which the initial resting state (i.e. the initial condition $\xi_0=0$, also the critical point of all maps in the quadratic family) remains bounded under iterations. For fixed $(c_0,c_1)$, the fixed-map Mandelbrot set designated the set of templates, visualized as a subset of the unit interval, for which the critical orbit stays bounded. These sets were found to be unions of subintervals in $[0,1]$, with structure tightly dependent on the fixed $(c_0,c_1)$ pair. We suggested that the distribution of these subintervals' lengths should be further investigated for power-law behavior. To study the Lebesgue measure of fixed-map Mandelbrot sets in $[0,1]$, we proceeded to define and analyze hybrid Mandelbrot sets.

For any fixed $c_0 \in \mathbb{C}$, the hybrid Mandelbrot slice was defined as an object in $\mathbb{C} \times [0,1]$, which can be visualized as a surface, or color plot: for any $c_1$ in the complex plane, the height/color is assigned based on the likelihood that a random template will keep the critical orbit bounded when iterated in conjunction with the pair of maps $(f_{c_0},f_{c_1})$. If, as proposed in the Introduction, we think of the fixed $f_{c_0}$ as the ``correct'' iteration map, the hybrid set ${\cal M}_{c_0}$ identifies, for each $c_1 \in \mathbb{C}$, how likely it is for the system initiated from rest to evolve along a sustainable (i.e., bounded) trajectory when randomly interposing ``erroneous'' maps $f_{c_1}$ in the iteration. In each $\mathbb{C} \times [0,1]$ hybrid plot, this likelihood can vary theoretically between 0 and 1. However, the maximum of 1 will not be necessarily achieved by all hybrid sets. As illustrated by Figure~\ref{hybrid_table}, the maxima of the hybrid sets generally get lower as $\lvert c_0 \rvert$ increases, but they do so in a non-trivial way.

To describe this phenomenon, we defined the contour Mandelbrot set, associating to each $c_0$ the highest probability (over all $c_1 \in \mathbb{C}$) for the template system to have bounded critical orbit. In particular, we considered the locus of $c_0 \in \mathbb{C}$ for which this highest probability value in the corresponding hybrid set attains the full value of 1. In other words, we searched for those ``correct'' maps $f_{c_0}$ for which there exist \emph{some} errors $f_{c_1}$ that deliver a sustainable evolution of the system initiated at rest, irrespective of the time of occurrence of these errors along the iteration. This locus contains the traditional Mandelbrot set, which warrants sustainability of the critical orbit for the system in absence of error.

We then studied how these results change when imposing the more restrictive condition that the orbits initiated at all critical points along the template iteration be bounded. We called the parameter loci obtained by this variation of the definition ``multi-critical Mandelbrot sets.'' As shown by standard theorems in non-autonomous dynamics, the multi-critical  Mandelbrot set for a fixed template iteration is equivalent with the Julia set connectedness locus. This equivalence has been extensively investigated by Hiroki Sumi, in the broader and somewhat different context of polynomial semigroups. In Sumi's framework, the answer is \emph{no}, and a counterexample is presented as Example 1.7 in~\cite{sumi2011dynamics}. However, the same reference provides a set of sufficient conditions for a postcritically bounded polynomial semigroup to have a connected Julia set (Theorem 2.14 in~\cite{sumi2011dynamics}).

We can formulate a practical interpretation of the equivalence of the Mandelbrot set with the Julia set connectedness locus. Suppose that the system is operating within the parameter range of the multi-critical Mandelbrot set; if the system is reset at any arbitrary iteration to its resting state (by which we mean the critical point $\xi_0=0$), this resting state is part of a connected prisoner set. In other words, the initial condition of the system can be perturbed continuously from rest into all other sustainable (i.e., asymptotically bounded) initial states.

Finally, we defined the multi-Mandelbrot set as the parameter locus in $\mathbb{C}^2$ of all $(c_0,c_1)$ which render the critical orbit bounded when iterated in any random order. This is locus has 4 real dimensions, making it hard to visualize per se. One way to represent it is to plot 2-dimensional slices obtained by fixing one of the two complex parameters (say, $c_0$), to obtain a collection of complex subsets corresponding precisely to the dark red hybrid plateaux illustrated in the panels of Figure~\ref{hybrid_table}. Another way to represent the multi-Mandelbrot set is via 3-dimensional slices obtained by fixing one parameter and some aspect of the second parameter, such as its modulus, or its real/imaginary part (as done in Figure~\ref{super_hybrid}). The multi-Mandelbrot set represents the set of all map pairs which produce a bounded critical orbit under all possible template iterations. In the context of applications, this represents all combinations of ``good'' and ``erroneous'' transformations which can be iterated in any template order while retaining a sustainable evolution for the critical set.

We noticed that the multi-Mandelbrot set is the same for both definitions. This is therefore a very desirable parameter range. First, if one thinks of $c_0$ as the correct map and of $c_1$ as its perturbation, or error, having a pair $(c_0,c_1)$ within the multi-Mandelbrot set guarantees sustainability of the resting system, independently of the frequency and timing of the errors along the template iteration. Second, this also implies sustainability if the system resets to rest at any point along the iteration.

\subsection{Future work}

This study represents only a first step in establishing and exploring significant questions for template iterations, in understanding their connection (similarities and differences) to the classical theory of single map iterations, as well as with exiting more general work in non-autonomous dynamics.

Our future work will focus on approaching analytically specific conjectures proposed in this  paper. For example, one aspect which we have not pursued in this paper is the change in the shape and boundary of hybrid Mandelbrot sets ${\cal M}_{c_0}$ as $c_0$ traverses different level sets of the contour set. Figure~\ref{hybrid_table} suggests that hybrid sets have increased fractal level sets when $0 < \beta(c_0) <1$ traverses the transitional contours of ${\cal M}$ (between the central plateau and the outside blue region). While this is intuitively not surprising, we would like to investigate the idea further, both computationally and analytically. A lot of theoretical work remains to be done, and exploring these questions may require a combination of methods from traditional and non-autonomous iterations.

Part of our future work is also aimed at understanding the significance of our results in the context of applications to the life sciences. In particular, we are considering the potential of using this theoretical framework to study natural iteration mechanisms, such as DNA replication. When a cell divides, it has to copy and transmit the exact same sequence of billions of nucleotides to its daughter cells. While most DNA is typically copied with high faithfulness, errors are a natural part of the process, and sometimes escape repair mechanisms, so that a mutated cell will end up being used as a template for the next replication iteration, with the possibility to lead to substantial, accumulated changes in the structure of later daughter cells. On one hand, accumulating mutations can lead to pathologies like cancer. On the other hand, perfect replication would lead to no genetic variation. Organisms may have to construct successful mechanisms that optimize between these two ends. In our previous work~\cite{radulescu2015symbolic}, we have suggested that template iterations may be appropriate to study how the size and timing of these mutations affect cells in the long term. In our current work, we are working to lay the grounds for contextualizing these questions within our mathematical framework, by introducing local errors in our template iterations (in the spirit of the errors made by chromosome replication mechanisms). 

\subsection*{Acknowledgements}

The work on this project was supported by SUNY New Paltz, via the Research Scholarship and Creative Activities Program, as well as the Research and Creative Projects Awards Program.

\bibliographystyle{unsrt}

\end{document}